\documentstyle[12pt]{article}
\begin{document}

\newtheorem{Def}{Definition}[section]
\newtheorem{thm}{Theorem}[section]
\newtheorem{lem}{Lemma}[section]
\newtheorem{rem}{Remark}[section]
\newtheorem{prop}{Proposition}[section]
\newtheorem{cor}{Corollary}[section]
\newtheorem{conj}{Conjecture}[section]
\newtheorem{question}{Question}[section]
\title
{ Local gradient estimates of solutions to some conformally
invariant fully nonlinear equations}
\author{YanYan Li\thanks{Partially
 supported by
        NSF grant DMS-0401118 and DMS-0701545.}\\
           Department of Mathematics\\
                  Rutgers University\\
                     110 Frelinghuysen Road\\
                        Piscataway, NJ 08854\\
                           USA
                              }
                          %\date{}
                          \maketitle
                          \input { amssym.def}

 \setcounter{section}{0}

\section{Introduction}

A classical theorem of Liouville says:
\begin{equation}
u\in C^2, \ \
\Delta u=0\ \mbox{and}\ u>0\ \mbox{in}\ \Bbb R^n
\ \mbox{imply that}\ u\equiv \ \mbox{constant}.
\label{1new}
\end{equation}

The Laplacian operator $\Delta$ is invariant under
rigid motions:  For any function u on $\Bbb R^n$ and
for any rigid motion $T:
\Bbb R^n\to \Bbb R^n$,
$$
\Delta (u\circ T)= (\Delta u)\circ T.
$$

$T$ is called a rigid motion if $Tx\equiv Ox+b$ for some $n\times n$
orthogonal matrix $O$ and some vector $b\in \Bbb R^n$.

It is clear that
 a linear second order partial differential operator
$$
Lu:= a_{ij}(x)u_{ij}+b_i(x)u_{i}+c(x)u
$$
is invariant under rigid motion, i.e.
$$
L(u\circ T)=(Lu)\circ T \ \mbox{for any function}\ u \ \mbox{and any
rigid motion}\ T,
$$
if and only if $L= a\Delta +c$ for some constants $a$ and $c$.

Instead of rigid motions, we look at M\"obius transformations of
$\Bbb R^n\cup \{\infty\}$ and {\it nonlinear} operators which are
invariant under M\"obius transformations.   A map $\varphi: \Bbb
R^n\cup \{\infty\}\to \Bbb R^n\cup \{\infty\}$ is called a M\"obius
transformation, if it is a composition of a finitely many of the
following three types of transformations:

\begin{eqnarray*}
 \mbox{A translation}: &
x\to x+\bar x,\  \mbox{where}\
\bar x
\ \mbox{ is a given point in}\   \Bbb R^n,
\\
\mbox{A dilation}:&
x\to a x, \ \mbox{where}\  a
\ \mbox{is a positive number},
\\
\mbox{ A Kelvin transformation}:&
x\to \frac x{|x|^2}.
\end{eqnarray*}

For a function $u$ on $\Bbb R^n$, let
$$
u_\varphi:= |J_\varphi|^{ \frac {n-2}{2n}} (u\circ\varphi)
$$
where
$J_\varphi$ denotes the Jacobian of $\varphi$.

Let $H(x, s, p, M)$ be a smooth function in its variables, where
$s>0$, $x, p\in \Bbb R^n$ and $M \in {\cal S}^{n\times n}$, the set
of $n\times n$ real symmetric matrices.  We say that a second order
fully nonlinear operator $H(\cdot, u, \nabla u, \nabla^2 u)$ is
conformally invariant if
$$
H(\cdot, u_\varphi, \nabla u_\varphi, \nabla^2 u_\varphi)\equiv
H(\cdot, u, \nabla u, \nabla^2 u)\circ \varphi
$$
holds for all positive smooth functions $u$ and all M\"obius
transformations $\varphi$.

For a positive $C^2$ function $u$, set
\begin{eqnarray}
A^u& := & -\frac{2}{n-2}u^{  -\frac {n+2}{n-2} } \nabla^2u+
\frac{2n}{(n-2)^2}u^ { -\frac {2n}{n-2} } \nabla u\otimes\nabla
u-\frac{2}{(n-2)^2} u^ { -\frac {2n}{n-2} }
|\nabla u|^2I\nonumber\\
&\equiv & A_w:=   w\nabla ^2 w-\frac {|\nabla w|^2}2 I, \qquad
\mbox{with}\ w=u^{ -\frac 2{n-2}},\label{AU}
\end{eqnarray}
where $I$ denotes the $n\times n$ identity matrix.

 Let $\varphi$ be a M\"obius
transformation, then for some $n\times n$ orthogonal matrix
functions $O(x)$ (i.e. $O(x) O(x)^t=I$), depending  on $\varphi$,
$$
A^{u_\varphi}(x)\equiv O(x)A^u(\varphi(x))O^t(x).
$$
Thus it is clear that  $f(\lambda(A^u))$ is a conformally invariant
operator for all symmetric functions $f$, where $\lambda(A^u)$
denotes the eigenvalues of $A^u$.

 It was proved in \cite{LL}
that an  operator $H(\cdot, u, \nabla u, \nabla^2 u)$ is conformally
invariant if and only if it is of the form
$$
H(\cdot, u, \nabla u, \nabla^2 u)\equiv f(\lambda(A^u)),
$$
where $f(\lambda)$ is some symmetric function in
$\lambda=(\lambda_1, \cdots, \lambda_n)$. Due to the above
characterizing conformal invariance property, the operator $A_w$ is
called the conformal Hessian of $w$.

Taking $f(\lambda)=\sigma_1(\lambda):=\lambda_1+\cdots+\lambda_n$,
we have                   a simple expression:
\begin{equation}
\sigma_1(\lambda(A^u))\equiv
-\frac 2{n-2}
 u^{-\frac{n+2}{n-2} }\Delta u.
\label{expression}
\end{equation}
In general, $f(\lambda(A^u))$ is
a  fully nonlinear operator, and is rather complex
even for $f(\lambda)=\sigma_k(\lambda)$, $k\ge 3$,  where
$$
\sigma_k(\lambda):=
\sum_{1\le i_1<\cdots<i_k\le n}\lambda_{i_1}\cdots \lambda_{i_k}
$$
is  the
$k-$th elementary symmetric function.  The expression for $\sigma_2$ is
still quite pleasant:
$$
\sigma_2(\lambda(A^u))\equiv \frac 12 \left(
\sigma_1(\lambda(A^u))^2 - (A^u)^t A^u \right).
$$

 Let
\begin{equation}
\Gamma\subset \Bbb R^n\
\mbox{be an open convex symmetric cone with vertex at the
origin}
\label{1-1}
\end{equation}
 satisfying
\begin{equation}
\Gamma_n:=\{\lambda\ |\ \lambda_i>0,
1\le i\le n\}\subset
\Gamma\subset \{\lambda\ |\ \sum_{i=1}^n\lambda_i>0\}=:\Gamma_1.
\label{1-2}
\end{equation}
Naturally, $\Gamma$ being symmetric means that
$(\lambda_1, \cdots, \lambda_n)\in \Gamma$
implies $(\lambda_{i_1}, \cdots, \lambda_{i_n})\in \Gamma$
for any permutation $(i_1, \cdots, i_n)$ of
$(1, \cdots,n)$.

 For
$1\le k\le n$,  let $\Gamma_k$ be
the connected component of $\{\lambda\in \Bbb R^n\ |\
\sigma_k(\lambda)>0\}$ containing the positive cone $\Gamma_n$.
It is known, see for instance \cite{CNS}, that $\Gamma_k$ satisfies
(\ref{1-1}) and (\ref{1-2}).   In fact
$$
\Gamma_k=\{\lambda\in \Bbb R^n\ |\ \sigma_1(\lambda), \cdots,
\sigma_k(\lambda)>0\}.
$$

Let $\Omega$ be an open subset of $\Bbb R^n$,  we consider
\begin{equation}\lambda(A^u)\in \partial \Gamma\qquad \mbox{in}\ \Omega,
\label{l1-1Prime} \end{equation} or
\begin{equation}
\lambda(A_w)\in \partial \Gamma\qquad \mbox{in}\ \Omega.
\label{L1-1}
\end{equation}
It is easy to see that in dimension $n\ge 3$
\begin{equation}
 A^u\equiv A_w\ \mbox{for any positive}\ C^2\
 \mbox{function}\ w\ \mbox{and}\
 u=w^{-\frac 2{n-2}}.
 \label{id}
\end{equation}

Equations (\ref{l1-1Prime}) and (\ref{L1-1})  are   fully nonlinear
second order degenerate elliptic equations. Fully nonlinear second
order elliptic equations with $\lambda(\nabla^2 u)$ in such general
$\Gamma$ were first studied by Caffarelli, Nirenberg and Spruck in
\cite{CNS}.

Equations (\ref{l1-1Prime}) and (\ref{L1-1}) have obvious meaning if
$u$ and $w$ are $C^2$ functions.   If they are in
$C^{1,1}_{loc}(\Omega)$, the equations are naturally understood to
be satisfied almost everywhere.   We give the
 notion of
viscosity solutions of (\ref{l1-1Prime}) and (\ref{L1-1}).

\begin{Def}
A positive continuous function $w$ in $\Omega$ is a viscosity
supersolution [resp. subsolution] of (\ref{L1-1}) when the following
holds: if $x_0\in \Omega$, $\varphi\in C^2(\Omega)$,
$(w-\varphi)(x_0)=0$ and $w-\varphi\ge 0$ near $x_0$ then
$$
\lambda(A_\varphi(x_0))\in \Bbb R^n\setminus \Gamma.
$$
[resp. if $(w-\varphi)(x_0)=0$ and $w-\varphi\le 0$ near $x_0$ then
$\lambda(A_\varphi(x_0))\in \overline \Gamma$].

We say that $w$ is a viscosity solution of (\ref{L1-1}) if it is
both a supersolution and a subsolution. \label{Def1.1}
\end{Def}

Similarly,  we have

\medskip

 \noindent {\bf
Definition \ref{Def1.1}$'$.}\ {\it A positive continuous function
$u$ in an open subset $\Omega$ of $\Bbb R^n$, $n\ge 3$, is a
viscosity subsolution [resp. supersolution] of
\begin{equation}\lambda(A^u)\in \partial \Gamma\qquad \mbox{in}\ \Omega,
\label{L1-1prime} \end{equation}
 when the following holds: if
$x_0\in \Omega$, $\varphi\in C^2(\Omega)$, $(u-\varphi)(x_0)=0$ and
$u-\varphi\le 0$ near $x_0$ then
$$
\lambda(A^\varphi(x_0))\in \Bbb R^n\setminus \Gamma.
$$
[resp. if $(u-\varphi)(x_0)=0$ and $u-\varphi\ge 0$ near $x_0$ then
$\lambda(A^\varphi(x_0))\in \overline \Gamma$].

We say that $u$ is a viscosity solution of (\ref{L1-1prime}) if it
is both a supersolution and a subsolution. }

\begin{rem}  In dimension $n\ge 3$,
a positive continuous function  $u$ is a viscosity subsolution
(supersolution) of (\ref{L1-1prime}) if and only if $w:=u^{ -\frac
2{n-2} }$ is a viscosity supersolution (subsolution) of
(\ref{L1-1}).  This is clear in view of (\ref{id}).
 \label{rem1.10}
\end{rem}

\begin{rem}  Viscosity solutions of (\ref{L1-1}) are invariant under conformal transformations
and  multiplication by positive constants.
 Namely, if $w$ is a viscosity supersolution (subsolution) of
(\ref{L1-1}) then, for any constants $b, \lambda>0$ and for any
$x\in \Bbb R^n$, $bw$ is a viscosity supersolution (subsolution) of
(\ref{L1-1}), $\xi(y):= \frac 1b w(x+b y)$ is a viscosity
supersolution (subsolution) of $\lambda(A_\xi)\in \partial \Gamma$
in $\{y\ |\ x+by \in \Omega\}$, and $\eta(y):= (\frac
{|y-x|^2}\lambda )^2 w(x+\frac{  \lambda^2(y-x) }{  |y-x|^2 })$ is a
viscosity supersolution (subsolution) of $\lambda(A_\eta)\in
\partial
\Gamma$ in $\{y\ |\ x+\frac{  \lambda^2(y-x) }{  |y-x|^2 }\in
\Omega\}$. \label{remA2}
\end{rem}

One of the two main theorems in this paper is the following Liouville theorem
for positive locally Lipschitz viscosity solutions of
\begin{equation}
\lambda(A_w)\in \partial \Gamma\qquad \mbox{in}\ \Bbb R^n.
\label{A3-1} \end{equation}

\begin{thm}
For $n\ge 3$, let $\Gamma$ satisfy (\ref{1-1}) and (\ref{1-2}), and
let $w$ be a positive locally Lipschitz viscosity solution of
(\ref{A3-1}).  Then $w\equiv w(0)$ in $\Bbb R^n$. \label{thm-vis}
\end{thm}

\begin{rem} For $n=2$, $\Gamma=\Gamma_1$, the conclusion
does not hold.  Indeed $w=e^{x_1}$ satisfies $\lambda(A_w)\in
\partial \Gamma_1$.  In fact, $\lambda(A_w)\in
\partial \Gamma_1$ is equivalent to $\Delta \log w=0$ in dimension $n=2$.
\label{rem1.0}
\end{rem}

Theorem \ref{thm-vis} can be viewed as a nonlinear extension of the
classical Liouville theorem (\ref{1new}). Indeed, in view of
(\ref{expression}),  Liouville theorem (\ref{1new}) is equivalent to
$$
u\in C^2, \ \ \lambda(A^u)\in \partial \Gamma_1 \ \mbox{and}\ u>0\
\mbox{in}\ \Bbb R^n \ \mbox{imply that}\ u\equiv \ \mbox{constant}.
$$

Such Liouville theorem was proved by  Chang, Gursky and Yang  in
\cite{CGY} for $u\in C^{1,1}_{loc}$,   $\Gamma=\Gamma_2$ and $n=4$;
by Aobing Li
 in \cite{Lia} for  $u\in C^{1,1}_{loc}$,   $\Gamma=\Gamma_2$ and $n=3$;
independently by Aobing Li
 in \cite{Lia} and by  Sheng, Trudinger and Wang in \cite{STW}
for  $u\in C^3$,   $\Gamma=\Gamma_k$, $k\le n$, $n\ge 3$. By
entirely different methods we established in \cite{Li2005c} the
following theorems.

Consider
\begin{equation}
f\in C^1(\Gamma)\cap C^0(\overline \Gamma)\ \mbox{is  symmetric in}\
\lambda_i, \label{1-3}
\end{equation}
\begin{equation}
f
 \ \mbox{is homogeneous of degree}\ 1,
\label{homo}
\end{equation}
\begin{equation}
f>0, \ f_{\lambda_i}:=\frac{\partial f}{\partial \lambda_i}>0 \
\mbox{in}\ \Gamma, \quad f|_{\partial \Gamma}=0, \label{1-4}
\end{equation}
\begin{equation}
\sum_{i=1}^n f_{\lambda_i} \ge \delta,\qquad \mbox{in}\ \Gamma\ \
\mbox{for some}\ \delta>0. \label{1-5new}
\end{equation}

Examples of such $(f, \Gamma)$ include  those given by elementary
symmetric functions:  For
 $1\le k\le n$, $(f, \Gamma)=(\sigma_k^{\frac 1k}, \Gamma_k)$
satisfies all the above properties; see for instance \cite{CNS}.
\medskip

\noindent{\bf Theorem A}\ ([24, v1])\ {\it For $n\ge 3$, let $(f,
\Gamma)$ satisfy (\ref{1-1}), (\ref{1-2}), (\ref{1-3}) and
(\ref{1-4}), and let $u$ be a positive $C^{1,1}_{loc}$ solution of
\begin{equation}
f(\lambda(A^u))=0,  \qquad \mbox{in}\ \Bbb R^n.
\label{liouville}
\end{equation}
Then $u\equiv u(0)$ in $\Bbb R^n$.
}

\medskip

\noindent{\bf Theorem B}\  ([24, final])\ {\it For $n\ge 3$, let
$(f, \Gamma)$ satisfy (\ref{1-1}), (\ref{1-2}), (\ref{1-3}) and
(\ref{1-4}), and let $u$ be a positive locally Lipschitz  weak
solution of (\ref{liouville}). Then $u\equiv u(0)$ in $\Bbb R^n$. }

\medskip

Throughout this paper, by a weak solution of (\ref{liouville}) we
mean in the sense of definition 1.1 in \cite{Li2005c}, with
$F(M):=f(\lambda(M))$ and $U:=\{M\ |\ \lambda(M)\in \Gamma\}$.   Our
proof of Theorem \ref{thm-vis} is along the line of
 \cite{Li2005c}, which makes use of ideas developed
 in \cite{LL2} and \cite{Li2005b} in treating the isolated singularity
 of $u$ at $\infty$.

\begin{rem} Let $(f, \Gamma)$ satisfy (\ref{1-1}),
(\ref{1-2}),
$$
f\in C^0(\overline \Gamma)\ \mbox{is symmetric in}\ \lambda_i,\ f>0\
\mbox{in}\ \Gamma, \ f|_\Gamma=0,
$$
$$
f(\lambda+\mu)\ge f(\lambda)\qquad \forall\ \lambda\in \Gamma,
\mu\in \Gamma_n,
$$
and let $\Omega$ be an open subset of $\Bbb R^n$. If $u$ is a
$C^{1,1}_{loc}$ solution of
\begin{equation}
f(\lambda(A^u))=0\qquad \mbox{in}\ \Omega, \label{K1}
\end{equation}
then it is a weak solution of (\ref{K1}). The proof is standard in
view of lemma 3.7 in \cite{Li2005c}. If $u$ is a weak solution of
(\ref{K1}), then it is clearly a viscosity solution of
 (\ref{L1-1prime}).
\end{rem}

The motivation of our study of such Liouville properties of
entire solutions of $\lambda(A^u) \in \partial \Gamma$
 is to answer the following questions
concerning local gradient estimates of solutions
to general second order conformally invariant
fully nonlinear elliptic equations.

Let $B_3\subset \Bbb R^n$ be a ball of radius $3$ and centered at
the origin.

\medskip

\noindent{\bf Question A}\ {\it Let $n\ge 3$, $(f, \Gamma)$ satisfy
(\ref{1-1}), (\ref{1-2}), (\ref{1-3})-(\ref{1-5new}). For constants
$0<b<\infty$ and
 $0<h\le  1$, let $u\in C^3(B_3)$ satisfy
 }
\begin{equation}
f(\lambda(A^u))=h, \ 0<u\le b, \ \lambda(A^u)\in \Gamma, \quad
\mbox{in}\ B_3. \label{A1}
\end{equation}
{\it Is it true that }
$$
|\nabla \log u|\le C\quad
\mbox{in}\ B_1 $$
{\it  for some constant $C$ depending only on  $b$  and $(f,
\Gamma)$? }

\medskip

Let $(M,g)$ be a smooth compact Riemannian manifold of dimension
$n\ge 3$. We use $i_0$ and  $R_{ijkl}$  to denote respectively the
injectivity radius and the curvature tensor. Consider the
 Schouten
 tensor
$$
A_g=\frac{1}{n-2}\left(Ric_g-\frac{R_g}{2(n-1)}g\right),
$$
where $Ric_g$ and $R_g$ denote respectively the
Ricci tensor and the scalar curvature.
  We use $\lambda(A_g)=(\lambda_1(A_g),
\cdots, \lambda_n(A_g))$ to denote the eigenvalues of
$A_g$ with respect to $g$.

Let $\hat g=u^{\frac 4 {n-2} }g$ be a conformal change of metrics,
then, see for example \cite{V00},
$$
A_{\hat g}=-\frac{2}{n-2}  u^{-1}\nabla^2u+ \frac{2n}{(n-2)^2} u^{-2}
\nabla u \otimes\nabla u -\frac{2}{(n-2)^2} u^{-2}
|\nabla u|^2g+A_{g},
$$
where covariant derivatives on the right-hand side are with respect to $g$.

For $g_1=u^{ \frac 4{n-2}}g_{flat}$, with $g_{flat}$ denoting the
Euclidean metric on $\Bbb R^n$,
$$
A_{g_1} =u^{ \frac 4{n-2} } A^u_{ij} dx^idx^j
$$
where $A^u$ is defined in (\ref{AU}).  In this case,
$\lambda(A_{g_1})=\lambda(A^u)$.

A more general question on Riemannian manifolds is

\medskip

\noindent{\bf Question B}\
{\it Let $g$ be a smooth Riemannian metric on $B_3\subset \Bbb R^n$,
$n\ge 3$, $(f, \Gamma)$ satisfy (\ref{1-1}), (\ref{1-2}),
(\ref{1-3})-(\ref{1-5new}). For a positive number $b$ and a positive
function $h\in C^1(B_3)$, let $u\in C^3(B_3)$ satisfy, with $\tilde
g:= u^{\frac 4{n-2}}g$,
\begin{equation}
f(\lambda(A_{\tilde g}))=h, \ 0<u\le b,\ \ \lambda(A_{\tilde g})\in
\Gamma, \quad \mbox{in}\ B_3. \label{abc1}
\end{equation}
Is it true that
\begin{equation}
\|\nabla \log u\|_g\le C  \quad \mbox{in}\ B_1
\label{abc2}
\end{equation}
for some constant $C$ depending only on $b, g$, $\|h\|_{ C^1(B_3) }$
and $(f, \Gamma)$? }

\medskip

 For $(f, \Gamma)=( \sigma_k^{\frac 1k}, \Gamma_k)$,
the local gradient estimate (\ref{abc2}) on Riemannian manifolds
 was established
by Guan and Wang in \cite{GW}; see a related work
 \cite{CGY} of Chang, Gursky and Yang where
 global
apriori $C^0$ and $C^1$  estimates for $f=\sigma_2^{\frac 12}$ and
$n=4$ were derived. Efforts  of achieving further generality were
made  in \cite{LL}, \cite{GW04}, 
\cite{STW} and \cite{GLW}.
 On locally
conformally flat manifolds, ``semi-local'' gradient estimates were
established, and used, in \cite{LL} and \cite{LL2} for  $(f,
\Gamma)$ satisfying (\ref{1-1}), (\ref{1-2}), (\ref{1-3}) and
(\ref{1-4}) via the method of moving spheres (or planes). A
consequence of the ``semi-local'' gradient estimates is, see also
lemma 0.5 and its proof in \cite{LL0},

\medskip

\noindent{\bf Theorem C}\ {\it Under an additional assumption $u\ge
a>0$ in $B_3$, the answer to Question A is ``Yes'', but with the
constant $C$ depending also on $a$. }

\begin{rem} In Theorem C, assumption (\ref{homo}) and (\ref{1-5new})
are not needed.
\end{rem}

Equations (\ref{A1}) and (\ref{abc1}) are fully nonlinear elliptic
equations of $u$.  Extensive studies have been given to fully
 nonlinear equations involving $f(\lambda(\nabla^2u))$ by
Caffarelli, Nirenberg and Spruck \cite{CNS}, Guan and Spruck
\cite{GS}, Trudinger \cite{T1}, Trudinger and Wang \cite{TW},  and
many others.

Fully nonlinear equations involving $f(\lambda(\nabla_g^2u+g))$ on
Riemannian manifolds are studied by Li \cite{Li90}, Urbas \cite{U},
and others.   Fully nonlinear equations on Riemannian manifolds
involving the Schouten tensor have been studied by Viaclovsky in
\cite{V3} and \cite{V1}, by Chang, Gursky and Yang in  \cite{CGY}
and \cite{CGY1}, and by many others; see for example \cite{CY},
\cite{Li1}, \cite{T}, \cite{V5}, and the references therein. Here we
study, on Riemannian manifolds $(M, g)$, local gradient estimates to
solutions of
\begin{equation}
f(\lambda(A_{ u^{  \frac 4{n-2} }g }))=h, \quad
\lambda(A_{ u^{  \frac 4{n-2} }g })\in \Gamma.
\label{local}
\end{equation}

If we make an additional concavity assumption
\begin{equation}
f\in C^2(\Gamma)\cap C^0(\overline \Gamma)\ \mbox{is symmetric in}\
\lambda_i, \mbox{and is
 concave in}\ \Gamma,
\label{s1}
\end{equation}
then we have the following corollary  of Theorem A and the proof of
(1.39) in \cite{LL}.

\begin{thm}  Let $(M,g)$ be as above and let $(f,\Gamma)$ satisfy
(\ref{1-1}), (\ref{1-2}), (\ref{homo}), (\ref{1-4}) and (\ref{s1}).
For a geodesic ball $B_{3r}$ in $M$ of radius $3r\le \frac 12 i_0$,
let $u$ be a $C^4$  positive solution
 of (\ref{local}) in $B_{3r}$.  Then
\begin{equation}
\|\nabla(\log u)\|_g\le C\qquad \mbox{in}\ B_r, \label{cc2}
\end{equation}
where $C$ is some positive constant depending only on $(f, \Gamma)$,
upper bounds of    $1/i_0$,  $\sup_{ B_{9r}}u$, $\|h\|_{ C^2(B_{9r})
}$  and a bound of $R_{ijkl}$ together with their covariant
derivatives up to second order. \label{thm2}
\end{thm}

 It has been observed
independently by Wang in \cite{W} that Theorem \ref{thm2} follows
from Theorem A.    The theorem is proved by Chen in \cite{C} using
a different method. It is well known, see e.g. \cite{CNS}, that $(f,
\Gamma)= (\sigma_k^{\frac 1k}, \Gamma_k)$ satisfies the hypotheses
of the theorem.

\begin{rem}  It is easy to see from Section 3 that
Theorem \ref{thm2} holds under slightly weaker hypotheses on $(f,
\Gamma)$:  Assuming that it satisfies (\ref{1-1}), (\ref{1-2}),
(\ref{1-4}), (\ref{s1}),
\begin{equation}
\lim_{s\to\infty} \ \inf_{\lambda\in K} f(s\lambda)=\infty\quad
\mbox{for any compact subset}\ K\ \mbox{of}\ \Gamma, \label{1-5}
\end{equation}
and
\begin{equation}
\inf_{  \lambda\in \Gamma, |\lambda|\ge \frac 1\delta } \left(
|\lambda|\sum_i f_{\lambda_i}(\lambda)\right) \ge \delta,\ \
\mbox{for some}\ \delta>0. \label{3-3}
\end{equation}
\label{aabbcc}
\end{rem}

The second main result in this paper is
\begin{thm} Let $(M,g)$ be as above and let $(f,\Gamma)$ satisfy
(\ref{1-1}), (\ref{1-2}), (\ref{1-3})-(\ref{1-5new}). For a geodesic
ball $B_{9r}$ in $M$ of radius $9r\le \frac 12 i_0$, let $u$ be a
$C^3$ positive solution
 of (\ref{local}) in $B_{9r}$.  Then (\ref{cc2}) holds,
where $C$ is some positive constant depending only on $(f, \Gamma)$,
upper bounds of   $1/i_0$, $\sup_{ B_{9r}}u$, $\|h\|_{ C^1(B_{9r})
}$
 and a bound
of $R_{ijkl}$ together with their first covariant derivatives.
\label{local2}
\end{thm}

\begin{rem} If $(f, \Gamma)$ satisfies (\ref{1-1}), (\ref{1-2}),
(\ref{1-3})-(\ref{1-4}), and $f$ is concave in $\Gamma$, then
(\ref{1-5new}) is automatically satisfied; see \cite{U}. Thus
Theorem \ref{local2} implies Theorem \ref{thm2}. The main point of
Theorem \ref{local2} is that no concavity assumption is made on $f$.
\end{rem}
\begin{rem}  Replacing the function $h$ in (\ref{local}) by $h(\cdot,
u)$ with $s^{\frac 4{n-2}}h(x, s)\in C^1(B_{9r}\times (0,
\infty))\cap L^\infty(B_{9r}\times (0, b))$ for all $b>1$, estimate
(\ref{cc2}) still holds, with the constant $C$ depending also on the
function $h$. This is easy to see form the proof of the theorem.
\label{moregeneral}
\end{rem}

\begin{rem}
Once (\ref{cc2}) is established,  it follows from
the proof of (1.39) in
 \cite{LL}, under the hypotheses of Theorem \ref{thm2},
 that
$$
\|\nabla_g^2(\log u)\|_g\le
C\qquad \mbox{in}\ B_r,
$$
where $C$ is some positive constant depending only on
 an
upper bound of   $1/i_0$,  $\sup_{ B_{9r}}u$,   $\sup_{ B_{3r}}
\|\nabla u\|_g$, $\|h\|_{ C^2(B_{9r}) }$  and a bound of $R_{ijkl}$
together with their covariant derivatives up to second order.
\label{rem0}
\end{rem}

A  subtlety of the local gradient estimate (\ref{cc2}) is that the
bound depends  on an upper bound of $u$, but not on  upper bounds of
$u^{-1}$. Global estimates of $|\nabla u|$ allowing the dependence
of an upper  bound of both $u$ and $u^{-1}$ was given by Viaclovsky
in \cite{V1}; see a related work \cite{Li90}. One application of the
local gradient estimate is for a rescaled sequence of solutions in
the following situation: For  solutions $\{u_i\}$ of (\ref{local})
in a unit ball $B_1$ satisfying, for some constant $b>0$ independent
of $i$,
$$
\sup_{B_1}u_i\le b u_i(0)\to \infty,
$$
consider
$$
v_i(y):= \frac 1{u_i(0)}v_i(\frac y{ u_i(0)^{ \frac 2{n-2} }}).
$$
One knows that
\begin{equation}
v_i(0)=1, \quad \mbox{and}
\ \ v_i(y)\le b\quad \forall\ |y|\le u_i(0)^{ \frac 2{n-2} },
\label{j1}
\end{equation}
and $v_i$ satisfies the same equation with $g$ replaced by the
rescaled metric $g^{(i)}$.
One would like to derive a bound of $|\nabla v_i|$
on $\{y\ |\  |y|<\beta\}$ for any fixed $\beta>1$.

Some time  ago the author arrived at the following idea:
Try to establish the estimate of  $|\nabla v_i|$
in two steps.

\medskip

\noindent{\it Step 1.}\
To establish,  for solutions $u$ of  (\ref{local}) for general
$(f, \Gamma)$,
local gradient estimates
 which depend on an upper
bound of both $u$ and $u^{-1}$.

\medskip

\noindent{\it Step 2.}\
To establish,  for solutions  $u$ of  (\ref{local})
in $B_1$ satisfying
$u(0)=1$,  an estimate on $B_\delta$
of  $u^{-1}$ from above,
  which depends on
 an upper bound of $u$.

\medskip

Once these two steps were achieved,
the needed gradient bound for solutions $\{v_i\}$ satisfying
(\ref{j1}) would follow.  The reason is that
we know from Step 2 that $v_i\ge a$ in
$B_\delta$ for some $a, \delta>0$ independent
of $i$.  Since $-L_{g^{(i)}}v_i\ge 0$
where $L_{g^{(i)}}$ denotes the conformal Laplacian
of $g^{(i)}$, and since
 $g^{(i)}$ tends to the Euclidean metric in
 $C^2_{loc}(\Bbb R^n)$, we have, for any $\beta>2$,
$$
v_i\ge  \xi_i\quad\mbox{on}\ B_\beta\setminus B_\delta,
$$
where $\xi_i$ is the solution of
$$
L_{g^{(i)}}\xi_i=0\ \ \mbox{in}\ B_\beta\setminus B_\delta,
\qquad \xi_i=  a\ \mbox{on}\ \partial B_\delta,\ \
\ \  \xi_i=  0 \ \mbox{on}\ \partial B_\beta.
$$
Clearly,
$$
\xi_i\to \frac {a\beta^{n-2}\delta^{n-2} }{ \beta^{n-2}-\delta^{n-2} }
\left( \frac 1{ |x|^{n-2} }- \frac 1 { \beta^{n-2} }\right)\quad
\mbox{uniformly in}\ B_\beta\setminus B_\delta.
$$
This  provides an upper bound of
$v_i^{-1}$ on $B_{\frac \beta 2}$, and the
desired estimate follows from Step 1.

Aobing Li and the author then started to implement this idea.  Step
1 for locally conformally flat manifolds was known to us, see
Theorem C. We established Step 1 on general manifolds and for
general $(f, \Gamma)$:

\medskip

\noindent{\bf Theorem D}\ {\it
\ (\cite{LL3})\
 Let $(M,g)$ be as above and let $(f,\Gamma)$ satisfy
(\ref{1-1}), (\ref{1-2}), (\ref{1-3})-(\ref{1-5new}). For a geodesic
ball $B_{9r}$ in $M$ of radius $9r\le \frac 12 i_0$, let $u$ be a
$C^3$ positive solution
 of (\ref{local}) in $B_{9r}$ satisfying,
for some positive constants $0<a<b<\infty$,
$$
a\le u\le b\qquad\mbox{on}\ B_{9r}.
$$
Then (\ref{cc2}) holds, where $C$ is some positive constant
depending only on $a$, $b$,  $\delta$, upper bounds of    $1/i_0$,
$\|h\|_{ C^1(B_{9r}) }$ and a bound of $R_{ijkl}$ together with
their first covariant derivatives. }

\medskip

This result was extended to manifolds with boundary under prescribed
mean curvature boundary conditions in \cite{JLL}; see theorem 1.3
there.  The proof of Theorem D uses Bernstein-type arguments.  The
choice of the auxiliary function $\phi$ in the proof is similar in
spirit to that in \cite{Li90} and \cite{V1}:  Finding a $\phi$ which
satisfies on a {\it finite interval} some second order ordinary
differential inequalities (see (\ref{claim})).  If the differential
inequalities (\ref{claim}) had a bounded solution $\phi$ on a half
line $(\alpha, \infty)$, then Theorem \ref{local2}, without the
assumption $u\ge a>0$,  would have been proved by the same method.
However the differential inequalities do not have any bounded
solution on any half line.

The method the author had in mind for Step 2 was to obtain, via
Bernstein-type arguments, a bound on
 $|\nabla \Phi(u)|=|\Phi'(u)\nabla u|$ for an appropriate $\Phi$.
For instance,
$|\nabla (u^\alpha)|\le C$ for $\alpha<0$ is weaker
than $|\nabla \log u|\le C$, and it becomes
weaker when $\alpha$ is smaller.
On the other hand, an estimate of $|\nabla (u^\alpha)|$ for any
$\alpha<0$ would yield an upper bound
of $u^{-1}$ near the origin.  In principal, estimating
$|\nabla (u^\alpha)|$ for very negative $\alpha$
should be easier than estimating $|\nabla \log u|$.
However we encountered some difficulties in
completing this step.

The author then took another path which requires establishing
appropriate Liouville theorems for  general degenerate conformally
invariant equations (\ref{liouville}). What needed is to prove that
any  positive  locally Lipschitz function $u$ satisfying
(\ref{liouville}) in appropriate weak sense must be a constant. In
[24, v1], a notion of weak solutions, tailored for the application
to local gradient estimates, was introduced.   Such Liouville
theorem for  $C^1_{loc}$ weak solutions of (\ref{liouville}) is
established there. My first impression was that weakening the
regularity assumption from $C^1_{loc}$ to $C^{0,1}_{loc}$ (locally
Lipschitz)
 is perhaps a subtle borderline issue whose solution
would  require some  new ideas
 beyond those
used in [24, v1]. It turns out,  to our surprise, that this only
requires some
 modification of our proof of the Liouville theorem
for $C^1_{loc}$ weak solutions. The  improvement, Theorem B, is
given  in [24, final].  Theorem B, together with Theorem C, is
enough to answer Question A affirmatively; this can be seen in the
proof of Theorem \ref{local2}.

With the help of the Jensen approximations (see \cite{J} and
\cite{CC}), we can further extend Theorem B for positive
locally Lipschitz viscosity solutions.  The theory of viscosity
solutions for nonlinear partial differential equations was developed
by Crandall and Lions in \cite{CL}.  Its basic idea also appears in
earlier papers by Evans \cite{E1, E2}.

Theorem \ref{thm-vis} allows us to,   using Theorem D, first
 establish a local  H\"older estimate
of $\log u$ instead of the local gradient estimate of $\log u$. With
the H\"older estimate of $\log u$, which yields the Harnack
inequality of $u$, we then obtain the local gradient estimate of
$\log u$ by another application of Theorem D.

The following problem looks
 reasonable and worthwhile to the author:
Using the Bernstein-type arguments to complete the above mentioned
Step 2, without any concavity assumption on $f$,  by choosing
appropriate $\Phi$.

One important ingredient   in our proof of Theorem \ref{thm-vis} is
 a new proof of the classical Liouville theorem
(\ref{1new}) which uses only the following two properties of
harmonic functions.

\medskip

\noindent{\it Conformal invariance of harmonic functions:}\ For any
harmonic function $u$, and for any M\"obius transformation
$\varphi$, $u_\varphi$ is harmonic.

\medskip

\noindent {\it  Comparison principle for harmonic functions on
balls:}\ Let $B\subset \Bbb R^n$, $n\ge 2$,  be the  ball centered
at the origin. Assume that $u\in C^2_{loc}(\overline
B\setminus\{0\})$ and $v\in C^2(\overline B)$ satisfy
$$
\Delta u= 0, \ \ u>0, \ \ \ \mbox{in}\ B\setminus\{0\},\qquad \qquad
\Delta v= 0\ \ \mbox{in}\ B,
$$
and
$$
u\ge v\qquad\mbox{on}\ \partial B.
$$
Then
$$
u\ge v \qquad \mbox{in}\ \overline B\setminus\{0\}.
$$

\bigskip

It is easy to see from this proof of the Liouville theorem
(\ref{1new}) that the
  following Comparison Principle is sufficient
  for a proof of Theorem \ref{thm-vis}.

\begin{prop} Let $\Omega\subset \Bbb R^n$, $n\ge 3$,  be a bounded open
set containing $m$ points $S_m:=\{P_1, \cdots, P_m\}$, $m\ge 0$,
$v\in C^{0,1}(\overline \Omega\setminus S_m)$, and $w\in
C^{0,1}(\overline \Omega)$.
 Assume that $w$ is a viscosity supersolution of
$ \lambda(A_w)\in \partial \Gamma$ in $\Omega\setminus S_m$, $v$ is
a viscosity subsolution of $ \lambda(A_v)\in \partial \Gamma$ in
$\Omega\setminus S_m$, and
$$
w>0\ \mbox{in}\ \overline \Omega, \ v>0\ \mbox{in}\ \overline
\Omega\setminus S_m, \  w>v\  \mbox{on}\
\partial \Omega.
$$
Then
\begin{equation}
\inf_{ \Omega \setminus S_m} (w-v)>0. \label{1c} \end{equation}
\label{propA5}
\end{prop}

\begin{rem} The proposition was proved in \cite{Li2005c}
under stronger hypotheses:  Instead of $C^{0,1}$ super or sub
viscosity solutions, they were assumed to be $C^{0,1}$ super or sub
weak solutions which include $C^{1,1}$ super or sub solutions.
\end{rem}

\begin{rem} Our equation $\lambda(A_w)\in \partial \Gamma$, or
(\ref{liouville}), does not satisfy the usual requirement on the
dependence on $w$ or $u$ in literature on viscosity solutions.
\end{rem}

\begin{rem}
The proof of Proposition \ref{propA5} for $m\ge 1$,
which makes use of the method developed in \cite{LL2}
(proof of theorem 1.3), \cite{Li2005b} (theorem 1.6-1.10)
and [24, final] (theorem 1.6 and remark 1.8) in treating
isolated singularities,  
is much more delicate than that for $m=0$, $S_0=\emptyset$.
 For $m=0$,  the
conclusion of the above theorem still holds in dimension $n=2$.  On
the other hand, the conclusion does not hold in dimension $n=2$ for
$\Gamma =\Gamma_1$ if $m\ge 1$. See the example below.
\label{rem1.1}
\end{rem}

\noindent {\it Example.}\ Let $ w(x)=(1+\epsilon)e^{-\frac 12x_1},
v(x)=e^{ -\frac 12 x_1 |x|^{-2} }$, $\epsilon>0$.  Clearly $w\in
C^\infty(\overline B_1)$, $v\in C^\infty(\overline
B_1\setminus\{0\})$, $w>v$ on $\partial B_1$, and they are positive
functions. Since $x_1$ is harmonic in $B_1$ and $x_1|x|^{-2}$ is
harmonic in $B_1\setminus \{0\}$, we know that $w\Delta w-|\nabla
w|^2=0$ in $B_1$ and $v\Delta v-|\nabla v|^2=0$ in $B_1\setminus
\{0\}$, i.e. $\lambda(A_w)\in \partial \Gamma_1$ in $B_1$ and $
\lambda(A_v)\in \partial \Gamma_1$ in $B_1\setminus\{0\}$. However,
 $\inf_{B_1\setminus\{0\}}(w-v)<0$ for small $\epsilon$.

\medskip

To prove Theorem \ref{thm-vis}, we only need Proposition
\ref{propA5} for $m=1$ and with $w\in C^{0,1}(\overline \Omega)$
being a viscosity supersolution of $\lambda(A_w)\in \partial\Gamma$
in $\Omega$.  In fact we only need a {\it weak comparison principle}
which assumes a priori $w\ge v$ in $\Omega\setminus\{0\}$; see
\cite{Li5}.

\begin{thm}
For $n\ge 3$, let $\Gamma$ satisfy (\ref{1-1}) and (\ref{1-2}), and
let $u$ be a positive locally Lipschitz viscosity solution of
\begin{equation}
\lambda(A^u)\in \partial \Gamma\qquad \mbox{in}\ \Bbb
R^n\setminus\{0\}. \label{AB} \end{equation}
 Then
\begin{equation}
u_{x,\lambda}(y) \le u(y),\qquad \forall\ 0<\lambda<|x|, |y-x|\ge
\lambda, y\ne 0. \label{JJJ}
\end{equation}
Consequently, $u$ is radially symmetric about the origin
 and
 $u'(r)\le 0$ for almost all $0<r<\infty$.
\label{thmA3new}
\end{thm}
The result was proved in \cite{Li2005c} under stronger hypotheses:
Assuming $u$ is a $C^{1,1}_{loc}$ or a $C^{0,1}_{loc}$ solution of
(\ref{AB}).

In the rest of the introduction we assume that $(M, g)$, $n\ge 3$,
is a smooth compact Riemannian manifold with nonempty smooth
boundary $\partial M$.  Let $h_g$ denote the mean curvature of
$\partial M$ with respect to the outer normal (a Euclidean ball has
positive mean curvature). For a conformal metric $\hat g=u^{\frac
4{n-2}}$, it is known that
$$
h_{\hat g}=u^{-\frac {n-2}2 }\left(-\frac{\partial u}{\partial
\nu_g} +\frac {n-2} 2 h_gu\right),
$$
where $\nu_g$ denotes the unit outer normal. We study
\begin{equation}
\left\{
\begin{array}{ll}
f(\lambda(A_{u^{\frac 4{n-2}}g}))=\psi, & \lambda(A_{u^{\frac
4{n-2}}g})\in \Gamma\ \mbox{on}\ O_1\setminus \partial M,\\
-\frac{ \partial u}{\partial \nu_g}+\frac {n-2}2 h_gu=\eta(x)u^{
\frac n{n-2}},& \qquad\quad \mbox{on}\ O_1\cap \partial M,
\end{array}
\right. \label{1.14}
\end{equation}
where $O_1$ is an open set of $M$, $\psi\in C^2(O_1)$ and $\eta\in
C^2(O_1\cap \partial M)$.

\begin{thm}\label{boundary} Assume that   $(M, g)$ is a smooth compact $n-$dimensional,
$n\ge 3$, Riemannian manifold with smooth boundary $\partial  M$,
and that  $(f,\Gamma)$ satisfy (\ref{1-1}), (\ref{1-2}),
(\ref{1-3})-(\ref{1-5new}). Let $O_1$ be an open set of $M$ and let
$u\in C^3(O_1)$ be a solution of (\ref{1.14}). If
$$
0< u\le b \quad \mbox{on } O_1 $$
 for some constant
$b$, then, for any open set ${O}_2$ of $ M$ satisfying
$\overline{O}_2\subset {O}_1$,
\begin{equation}
|\nabla (\log u)|_g\le C \quad \mbox{on } {O}_2 \label{cc9}
\end{equation}
 for some positive constant $C$ depending only on $n$ $(f,
\Gamma)$, $(M, g)$, $\psi$, $\eta$, $b$, $O_1$ and $O_2$.
\end{thm}

\begin{rem}  When $(f, \Gamma)$ satisfies a more restrictive condition $(H_1)$
 defined  in \cite{LL}, which includes all $(\sigma_k^{1/k}, \Gamma_k)$,
estimate (\ref{cc9}) was established in \cite{JLL}.
\end{rem}

\begin{rem}  Replacing the function $\psi$ and $\eta$  in (\ref{1.14})
respectively by $\psi(\cdot, u)$ and $\eta(\cdot, u)$ satisfying
$s^{\frac 4{n-2}}\psi(x,s)\in C^2(O_1\times (0, \infty)) \cap
L^\infty(O_1\times (0, b)))$ and $\eta\in C^2((O_1\cap
\partial M)\times (0, \infty))\cap
L^\infty((O_1\cap
\partial M)\times (0, b))$ for all
$b>1$,  estimate (\ref{cc9}) still holds. This is easy to see from
the proof of the theorem. \label{moregeneral1}
\end{rem}

Let
$$
\Bbb R^n_+:=\{ x=\{x', x_n)\in \Bbb R^n\ |\ x'=(x_1, \cdots,
x_{n-1}), x_n>0\}
$$
denote the half Euclidean space, and let $\Omega^+\subset \Bbb
R^n_+$ be an open set.  We use notations
$$
\partial''\Omega^+=\overline{ \partial \Omega^+\cap \Bbb R^n_+ },
\quad \partial'\Omega^+=\partial \Omega^+\setminus \partial''
\Omega^+.
$$

The following definition is standard.
\begin{Def}
A function $u\in C^0(\overline {\Omega^+})$ is said to satisfy
$$
\frac{\partial u}{\partial x_n}\le 0 \ \ [\mbox{resp.}\ \ge 0]
\qquad\mbox{on}\ \partial'\Omega^+
$$
in the viscosity sense, if $\bar x\in\partial'\Omega^+$, $\psi\in
C^1(\overline {\Omega^+})$ and $u-\psi$ has a local minimum [resp.
local maximum] at $\bar x$ then
\begin{equation}
\frac {\partial \psi}{ \partial x_n}(\bar x)\le 0 \ [\mbox{resp.}\
\ge 0]. \label{13a}
\end{equation}
Similarly we define
$$
\frac{\partial u}{\partial x_n}<0\ \ \mbox{or}\ \ \frac{\partial
u}{\partial x_n}>0 \ \ \mbox{on}\
\partial'\Omega^+\
\mbox{in the viscosity sense}
$$
by making the inequalities in (\ref{13a}) strict.

We say that $\frac{\partial u}{\partial x_n}=0$ on
$\partial'\Omega^+$ in the viscosity sense if both $\frac{\partial
u}{\partial x_n}\le 0$ and $\frac{\partial u}{\partial x_n}\ge 0$ on
$\partial'\Omega^+$ in the viscosity sense.
\end{Def}

\begin{thm}
Let $\Gamma$ satisfy (\ref{1-1}) and (\ref{1-2}), and let $u\in
C^{0,1}(\overline{\Bbb R^n_+})$ be a positive viscosity solution of
\begin{equation}
\lambda(A^u)\in \partial \Gamma\qquad \mbox{in}\ \Bbb R^n_+
\label{16a}
\end{equation}
satisfying, in the viscosity sense \begin{equation} \frac{\partial
u}{\partial x_n}=0 \qquad \mbox{on}\ \partial \Bbb R^n_+.
\label{16b} \end{equation} Then $u\equiv u(0)$ in $\overline{\Bbb
R^n_+}$. \label{thm16}
\end{thm}

Theorem \ref{thm-vis}, Theorem \ref{thm2} and Theorem \ref{local2}
were announced in \cite{Li4}, and the proofs were given in
\cite{Li2006a}.  The proof of Theorem \ref{thm-vis} in this revised
version of \cite{Li2006a} is improved in presentation.

 The paper is organized as follows. In Section 2 we prove
Proposition \ref{propA5}, Theorem \ref{thm-vis} and Theorem
\ref{thmA3new}. In Section 3 we prove Theorem \ref{thm2}. In Section
4 we prove Theorem \ref{local2}. In Section 5 we prove Theorem
\ref{thm16} and Theorem \ref{boundary}. In Appendix A we give, for
reader$'$s convenience, the proof of Theorem D.

We end the introduction by a question related to Theorem
\ref{thm-vis}.   Let  $\Gamma$ satisfy (\ref{1-1}) and (\ref{1-2}),
and let  $E\in C^\infty(\Bbb R^n\times \Bbb R_+\times \Bbb R^n)$
satisfy
$$
E(x,\alpha s,\alpha p)= \alpha E(x, s, p)\quad \forall \ \alpha,
s>0, x, p\in \Bbb R^n,
$$
and
$$
E(x, 1, 0)\equiv 0\quad \forall\ x\in \Bbb R^n.
$$
Assume that $w$ is a positive function in $C^\infty(\Bbb R^n)$
satisfying
\begin{equation}
\lambda(\nabla ^2 w+E(x, w, \nabla w))\in \partial \Gamma,\qquad
\forall\ x\in \Bbb R^n. \label{nn7}
\end{equation}

\begin{question} Under what additional hypothesis on $E$, the above
imply that $w\equiv w(0)$ on $\Bbb R^n$?  What if $w$ has less
regularity, e.g. in $C^{1,1}_{loc}(\Bbb R^n)$, or a locally
Lipschitz viscosity solution of (\ref{nn7})?
\end{question}

We know from Theorem \ref{thm-vis} and Remark \ref{rem1.0} that for
$E(x, w, \nabla w)\equiv -\frac{|\nabla w|^2}{2w}I$ the answer is
``yes'' in dimension $n\ge 3$ and is ``No'' in dimension $n=2$. What
about $E(x, w, \nabla w)\equiv b\frac {|\nabla w|^2}wI$ for other
constants $b$?

\section{Proof of Proposition
\ref{propA5}, Theorem \ref{thm-vis} and Theorem \ref{thmA3new}}

We first give

\noindent{\bf A new proof of the classical Liouville theorem
(\ref{1new})}.\ For every $x\in \Bbb R^n$, and for every
$\lambda>0$, let
$$
u_{x, \lambda}(y):= \frac{ \lambda^{n-2}}{  |y-x|^{n-2} } u(x+\frac{
\lambda^2(y-x) }{  |y-x|^2  }).
$$
We know that $u_{x, \lambda}=u$ on $\partial B_\lambda(x)$, $u\in
C^0(B_\lambda(x))$, $u_{x, \lambda}\in
C^0(B_\lambda(x)\setminus\{x\})$, $u$ and $u_{x, \lambda}>0$ are
positive harmonic functions in
 $B_\lambda(x)$ and  $B_{\lambda}(x)\setminus\{x\}$ respectively. Note that
 we have used the conformal invariance of harmonic functions to obtain
 the harmonicity of $u_{x,\lambda}$.  By
 the comparison principle for harmonic functions
 on balls, $u_{x, \lambda}\ge u$ in $B_\lambda(x)\setminus\{x\}$
 which is equivalent to $u_{x,\lambda}\le u$ in
 $\Bbb R^n\setminus B_\lambda(x)$. It
 follows that $u\equiv u(0)$; see e.g. lemma 11.2 in \cite{LZ} or
 lemma A.1 in \cite{LL2}.

Now we give

\noindent{\bf The proof of Theorem \ref{thm-vis} using Proposition
\ref{propA5}.}\  For every $x\in \Bbb R^n$, and for every
$\lambda>0$, applying Proposition \ref{propA5} to $u$ and $u_{x,
\lambda}$ on $B_{\lambda}(x)$ yields $u_{x, \lambda}\ge u$ in
$B_\lambda(x)\setminus\{x\}$. This implies $u\equiv u(0)$.

Here is

\noindent{\bf The proof of Theorem \ref{thmA3new}
 using Proposition
\ref{propA5}.}\  For every $x\in \Bbb R^n\setminus \{0\}$, and for
every $0<\lambda<|x|$, applying Proposition \ref{propA5} to $u$ and
$u_{x, \lambda}$ on $B_{x,\lambda}(x)$ yields $u_{x, \lambda}\ge u$
in $B_\lambda(x)\setminus\{x, |x|^{-2}(|x|^2-\lambda^2)x\}$, i.e.
(\ref{JJJ}) holds. It follows that $u$ is radially symmetric about
the origin
 and
 $u'(r)\le 0$ for almost all $0<r<\infty$; see e.g. \cite{Li2005b}.

\medskip

In the rest of this section we give the

\noindent{\bf Proof of   Proposition \ref{propA5}}.\
 We prove it by induction on
the number of points $m$.  We start from $m=0$ with $S_0=\emptyset$.

\medskip.

 \noindent{\it Step 1}. \  Proposition \ref{propA5} holds
for $m=0$.

\medskip

 Because of Remark \ref{remA2}, we only need to show that
$w\ge v$ in $\Omega$. We prove it by contradiction. Suppose the
contrary, then, for some $\gamma>0$,
$$
\max_\Omega(v-w)\ge \gamma, \quad v-w\le -\gamma\ \mbox{on}\
\Omega\setminus \Omega_\gamma,
$$
where $\Omega_\gamma:= \{ x\in \Omega\ |\ dist(x, \partial \Omega)
>\gamma\}$.

For small positive constants $0<\epsilon<<\mu<<\delta<<1$ which we
specify later, let
\begin{equation}
\hat v(y):= v(y)+\mu \varphi(y), \qquad \varphi(y):= e^{ \delta
|y|^2}, \label{B2-1}
\end{equation}
$$
\hat v^\epsilon (y):= \sup_{ x\in \Omega} \{ \hat v(x) +
\epsilon-\frac 1\epsilon |x-y|^2\}, $$
$$
 w_\epsilon(y):= \inf_{ x\in \Omega} \{ w(x) -
\epsilon+\frac 1\epsilon |x-y|^2\}. $$

$\hat v$ has been used in \cite{Li2005c}.  $\hat v^\epsilon$ and
$w_\epsilon$ are Jensen approximations whose useful properties can
be found in theorem 1.5 and lemma 5.2 in \cite{CC}.  We list below
some properties which we need.
\begin{equation}
\hat v^\epsilon\ \mbox{and}\ w_\epsilon\ \mbox{are punctually second
order differentiable a.e. in }\ \Omega_\gamma, \label{B3-2}
\end{equation}
\begin{equation}
\nabla^2 \hat v^\epsilon\ge -\frac C\epsilon I, \ \ \nabla^2
w_\epsilon \le \frac C\epsilon I, \qquad \mbox{in}\ \Omega_\gamma,
\label{B3-1}
\end{equation}

For any $x\in \Omega_\gamma$, there exist $x^*=x^*(x)$ and
$x_*=x_*(x)$ in $\Omega$ such that
\begin{equation}
\hat v^\epsilon(x)=\hat v(x^*) +\epsilon-\frac 1\epsilon |x^*-x|^2,
\label{B4-1}
\end{equation}
\begin{equation}
w_\epsilon(x)=w(x_*)-\epsilon+\frac 1\epsilon |x_*-x|^2,
\label{B4-2}
\end{equation}
\begin{equation}
|x^*-x| +|x_*-x|\le C\epsilon, \label{B4-3}
\end{equation}
\begin{equation}
|\nabla \hat v^\epsilon|+|\nabla w_\epsilon|\le C \quad \mbox{in}\
\Omega_\gamma, \label{B4-4}
\end{equation}
where, and in the following, $C$ denotes various positive constants
independent of $\mu, \delta$ and $\epsilon$. The punctual second
order differentiability is defined as in definition 1.4 in
\cite{CC}. Properties (\ref{B3-2})-(\ref{B4-2}) can be found in
\cite{CC} which hold for continuous $w$ and $v$. Property
(\ref{B4-3}) follows from the proof of (5) in lemma 5.2 in \cite{CC}
by using the Lipschitz regularity of $w$ and $v$. Property
(\ref{B4-4}) can easily be deduced as follows from
(\ref{B4-1})-(\ref{B4-3}) using again the Lipschitz regularity of
$w$ and $v$: For any $y, z$ in $\Omega_\gamma$, we have, by
(\ref{B4-1}) and the definition of $\hat v^\epsilon$,
\begin{eqnarray*}
\hat v^\epsilon(y)&\ge & \hat v(y-z+z^*)+\epsilon-\frac 1\epsilon
|z^*-z|^2 =\hat v(y-z+z^*)-\hat v(z^*)+\hat v^\epsilon(z)\\
&\ge & \hat v^\epsilon(z)-C|y-z|.
\end{eqnarray*}
This gives the bound of $|\nabla \hat v^\epsilon|$ in (\ref{B4-4}).
The bound of $|\nabla w_\epsilon|$ can be obtained similarly.

Using the Lipschitz regularity of $v$ and $w$, it is easy to deduce
from (\ref{B4-1})-(\ref{B4-3}) that
$$
|\hat v^\epsilon-\hat v|+|w_\epsilon-w|\le C\epsilon \quad
\mbox{in}\ \Omega_\gamma.
$$
Thus, for small $\epsilon$ and $\mu$, there exists $1<b_\epsilon<C$
such that
$$
\hat v^\epsilon-b_\epsilon w_\epsilon \le -\gamma\quad \mbox{in}\
\Omega\setminus \Omega_\gamma,
$$
$$
\max_\Omega\{\hat v^\epsilon-b_\epsilon w_\epsilon\}=\epsilon.
$$
Let $\xi_\epsilon:= \hat v^\epsilon-b_\epsilon w_\epsilon$, and let
$\Gamma_{  \xi_\epsilon^+ }$ denote the concave envelope of
$\xi_\epsilon^+:= \max\{ \xi_\epsilon, 0\}$ on $\Omega$. By
(\ref{B3-1}),
$$
\nabla^2 \xi_\epsilon\ge -\frac C\epsilon I\qquad \mbox{in }\
\Omega_\gamma.
$$
Thus, by lemma 3.5 of \cite{CC},
$$
\int_{  \{\xi_\epsilon=\Gamma_{ \xi_\epsilon^+}  \}  } \det
(-\nabla^2 \Gamma_{ \xi_\epsilon^+})>0.
$$
It follows that the Lebesgue measure of $\{\xi_\epsilon=\Gamma_{
\xi_\epsilon^+}  \}$ is positive.  By (\ref{B3-2}), there exists
$x_\epsilon\in \{\xi_\epsilon=\Gamma_{ \xi_\epsilon^+}  \}$ such
that  both $\hat v^\epsilon$ and $w_\epsilon$ are punctually second
order differentiable at $x_\epsilon$.  Clearly, for small
$\epsilon$, $x_\epsilon\in \Omega_\gamma$,
\begin{equation}
0<\xi_\epsilon(x_\epsilon)<\epsilon, \label{B8-1}
\end{equation}
\begin{equation}
|\nabla \xi_\epsilon(x_\epsilon)|\le C\epsilon, \label{B8-2}
\end{equation}
\begin{equation}
\nabla ^2 \xi_\epsilon(x_\epsilon)= \nabla ^2 \hat
v^\epsilon(x_\epsilon)-b_\epsilon \nabla ^2 w_\epsilon(x_\epsilon)
\le 0, \label {B8-3}
\end{equation}
\begin{equation}
w_\epsilon(x_\epsilon+z)\ge w_\epsilon(x_\epsilon) +\nabla
w_\epsilon(x_\epsilon) \cdot z +\frac 12 z^t \nabla ^2
w_\epsilon(x_\epsilon)z -\circ(|z|^2), \label{B9-1}
\end{equation}
\begin{equation}
\hat v^\epsilon(x_\epsilon+z)\le \hat v^\epsilon(x_\epsilon) +\nabla
\hat v^\epsilon(x_\epsilon) \cdot z +\frac 12 z^t \nabla ^2 \hat
v^\epsilon(x_\epsilon)z +\circ(|z|^2), \label{B9-2}
\end{equation}
By the definition of $\hat v^\epsilon$, we have, with
$(x_\epsilon)^* =(x_\epsilon)^*(x)$ as in (\ref{B4-1}),
$$
\hat v^\epsilon(x_\epsilon+z)\ge \hat v( (x_\epsilon)^*+z)
 +\epsilon -\frac 1\epsilon |(x_\epsilon)^*-x_\epsilon|^2
$$
and therefore, in view of (\ref{B9-2}) and (\ref{B2-1}),
\begin{eqnarray*}
v( (x_\epsilon)^*+z)&\le &\hat v^\epsilon(x_\epsilon+z) -\epsilon+
\frac 1\epsilon |(x_\epsilon)^*-x_\epsilon|^2-\mu \varphi(
(x_\epsilon)^*+z)\\
&\le & Q_\epsilon(z)+\circ(|z|^2),
\end{eqnarray*}
where $Q_\epsilon(z)$ is the quadratic polynomial with
$$
Q_\epsilon(0)=\hat v^\epsilon(x_\epsilon)-\epsilon +\frac 1\epsilon
|(x_\epsilon)^*-x_\epsilon|^2- \mu \varphi( (x_\epsilon)^*) = \hat
v^\epsilon(x_\epsilon) - \mu \varphi( (x_\epsilon)^*)+O(\epsilon),
$$
$$
\nabla Q_\epsilon(0)=\nabla \hat v^\epsilon(x_\epsilon)- \mu\nabla
\varphi( (x_\epsilon)^*),
$$
$$
\nabla ^2 Q_\epsilon(0)=\nabla^2 \hat v^\epsilon(x_\epsilon)-
\mu\nabla ^2\varphi( (x_\epsilon)^*),
$$
where $|O(\epsilon)|\le C\epsilon$.

By (\ref{B4-1}) and (\ref{B2-1}), $Q_\epsilon(0)=v( (x_\epsilon)^*
)$.  Since $v$ is a viscosity subsolution of (\ref{L1-1}), we have
\begin{equation}
\lambda (A_{Q_\epsilon}(0))\in \overline \Gamma. \label{B10-1}
\end{equation}

For small $0<\epsilon<< \mu<<\delta$, we have, as in the proof of
lemma 3.7 in \cite{Li2005c}, that
\begin{equation}
A_{Q_\epsilon}(0)\le  ( 1-\mu \frac {  [\varphi(
(x_\epsilon)^*)+O(\frac \epsilon\mu)] } {  \hat
v^\epsilon(x_\epsilon) })A_{  \hat v^\epsilon}(x_\epsilon)  -\frac
{\mu\delta}2   [\varphi( (x_\epsilon)^*)+O(\frac \epsilon{\mu
\delta}+\frac \mu\delta)] \hat
v^\epsilon(x_\epsilon) I. \label{B11-1}
\end{equation}
Similarly, using (\ref{B9-1}) and the definition of $w_\epsilon$, we
have
$$
w( (x_\epsilon)_*+z)\ge w_\epsilon(x_\epsilon+z)+\epsilon-\frac
1\epsilon |(x_\epsilon)_*-x_\epsilon|^2 \ge P_\epsilon(z)
-\circ(|z|^2),
$$
where $P_\epsilon(z)$ is the quadratic polynomial with
$$
P_\epsilon(0)=w_\epsilon(x_\epsilon)+\epsilon - \frac 1\epsilon
|(x_\epsilon)_\epsilon-x_\epsilon|^2=
w_\epsilon(x_\epsilon)+O(\epsilon),
$$
$$\nabla P_\epsilon(0)=\nabla w_\epsilon(x_\epsilon),
\quad \nabla^2 P_\epsilon(0)=\nabla ^2 w_\epsilon(x_\epsilon).
$$
By (\ref{B4-2}), $P_\epsilon(0)= w(  (x_\epsilon)_*)$. Since $w$ is
a viscosity supersolution of (\ref{L1-1}), we have
\begin{equation}
\lambda(A_{ P_\epsilon }(0))\in \Bbb R^n\setminus \Gamma.
\label{B12-1}
\end{equation}
By (\ref{B4-4}),
\begin{equation}
A_{ P_\epsilon(0) }= [w_\epsilon(x_\epsilon)+O(\epsilon)] \nabla^2
w_\epsilon(x_\epsilon) -\frac 12 |\nabla w_\epsilon(x_\epsilon)|^2
I= (1+O(\epsilon)) A_{ w_\epsilon}(x_\epsilon)+O(\epsilon).
\label{B13-1}
\end{equation}
By (\ref{B8-3}),
$$
A_{ w_\epsilon}(x_\epsilon)\ge \frac {  w_\epsilon(x_\epsilon) } {
b_\epsilon  \hat v^\epsilon(x_\epsilon)  } A_{  \hat v^\epsilon
}(x_\epsilon) + \frac 1{  2 b_\epsilon^2 } (  \frac {   b_\epsilon
w_\epsilon(x_\epsilon)  } {\hat v^\epsilon (x_\epsilon)  } |\nabla
\hat v^\epsilon(x_\epsilon)|^2 - \nabla (b_\epsilon w_\epsilon)
(x_\epsilon)|^2 ) I.
$$
By (\ref{B8-1}) and (\ref{B8-2}),
$$
|\frac {   b_\epsilon w_\epsilon(x_\epsilon)  } {\hat v^\epsilon
(x_\epsilon)  }   -1|\le C\epsilon, \quad |\nabla \hat
v^\epsilon(x_\epsilon)-
 \nabla (b_\epsilon w_\epsilon)
(x_\epsilon)|\le C\epsilon.
$$
It follows, in view of (\ref{B4-4}),
$$
| \frac {   b_\epsilon w_\epsilon(x_\epsilon)  } {\hat v^\epsilon
(x_\epsilon)  } |\nabla \hat v^\epsilon(x_\epsilon)|^2 - |\nabla
(b_\epsilon w_\epsilon) (x_\epsilon)|^2|\le C\epsilon,
$$
and
\begin{equation}
A_{ w_\epsilon}(x_\epsilon)\ge \frac {  w_\epsilon(x_\epsilon) } {
b_\epsilon  \hat v^\epsilon(x_\epsilon)  } A_{  \hat v^\epsilon
}(x_\epsilon)  -C\epsilon I. \label{B14-1}
\end{equation}

By (\ref{B13-1}), (\ref{B14-1}) and (\ref{B11-1}), we have, after
fixing some small $0<\mu<<\delta$,
$$
A_{P_\epsilon}(0)\ge a(\mu, \delta) A_Q(0)+b(\mu, \delta)I-C(\mu,
\delta)\epsilon I,
$$
where $a(\mu, \delta)$, $b(\mu, \delta)$ and $C(\mu, \delta)$ are
some positive constants independent of $\epsilon$.  Now fix
$\epsilon>0$ such that $b(\mu, \delta)-C(\mu, \delta)\epsilon>0$, we
deduce from (\ref{B10-1}), using the properties of $\Gamma$,  that
$\lambda(A_{P_\epsilon}(0))\in\Gamma$.  This violates (\ref{B12-1}).
Step 1 is established.

\medskip

\noindent{\it Step 2.}\ Proposition \ref{propA5} holds for $m$ if it
holds for $m-1$.

\medskip

 Now we assume that the proposition holds for $m-1$
points, $m-1\ge 0$, and we will prove that it holds for $m$ points.
We prove (\ref{1c}) by contradiction.  Suppose it does not hold,
then
$$
\inf_{\Omega\setminus S_m}(w-v)\le 0.
$$
By Shrinking $\Omega$ slightly, and working with the smaller one, we
may assume without loss of generality that $w$ is $C^{0,1}$ in some
open neighborhood of $\overline\Omega$.   Let
$$
u:= v^{ -\frac {n-2}2 }\qquad \mbox{and}\qquad \xi:= w^{ -\frac
{n-2} 2}.
$$
Then $$
\inf _{ \Omega\setminus S_m}(u-\xi)\le 0,
\qquad u>\xi\ \ \mbox{on}\
\partial\Omega,
$$
$u$ is a viscosity supersolution of
\begin{equation}
\lambda(A^u)\in \partial \Gamma\quad\mbox{in}\ \Omega\setminus S_m,
\label{3b}
\end{equation}
and $\xi$ is a viscosity subsolution of $$
\lambda(A^\xi)\in \partial\Gamma\qquad \mbox{in}\ \Omega. $$

For a positive $C^2$ function $\psi$, $A^\psi(x_0)\in
\overline\Gamma$ implies $\Delta\psi(x_0)\le 0$.  So by the
definition of $u$ being a viscosity supersolution of (\ref{3b}),
$$
\Delta u\le 0\ \mbox{in}\ \Omega\setminus S_m\ \mbox{in the
viscosity sense}.
$$
It follows, using also the positivity
of $u$,  that
$$
\inf_{ \Omega\setminus S_m} u\ge \inf_{\partial \Omega}u>0.
$$
Thus, for some $0<a\le 1$,
$$
\inf_{\Omega\setminus S_m} (u-a\xi)=0.
$$
Since we can use $a^{-1}u$ instead of $u$, we may assume without
loss of generality that $a=1$.  So we have, in addition,
$$
\inf_{\Omega\setminus S_m} (u-\xi)=0. $$
Let $P_m$ be the origin, and let
$$
\widehat \Omega:= \Omega\setminus S_{m-1}, \ S_{m-1}:=\{P_1, \cdots,
P_{m-1}\}.
$$
\begin{lem}
There exists $\epsilon>0$ such that $u=\xi=\xi(0)$ in
$B_\epsilon\setminus\{0\}$. \label{lemA5}\end{lem}

\noindent{\bf Proof of Lemma \ref{lemA5}.}\ We first claim that
\begin{equation}
\liminf_{ |y|\to 0}(u-\xi)(y)=0. \label{5b}
\end{equation}
Indeed if (\ref{5b}) did not hold, there would be some $\epsilon>0$
such that $\displaystyle{ \inf_{B_\epsilon\setminus\{0\}}(u-\xi)>0
}$ i.e.
\begin{equation}
\inf_{B_\epsilon\setminus\{0\}} (w-v)>0. \label{5c}
\end{equation}
Since the singular set of $v$ in
$\Omega\setminus\overline{B_\epsilon}$ is $S_{m-1}$ which contains
only $m-1$ points, we have, by the induction hypothesis,
$$
\inf_{    (\Omega\setminus B_\epsilon)\setminus S_{m-1}  } (w-v)>0.
$$
This and (\ref{5c}) violate (\ref{5b}).

Now let
$$
\Phi(\xi, x, \lambda; y):=\lambda \xi(x+y).
$$
Since $\Phi(\xi, 0, 1; \cdot)=\xi$ and $u>\xi$ on $\partial \Omega$,
we can fix some $\epsilon_4>0$ so that $|x|\le \epsilon_4$ and
$|\lambda-1|\le \epsilon_4$ guarantee
\begin{equation}
u>\Phi(\xi, x, \lambda; \cdot)\qquad \mbox{on}\ \partial \Omega.
\label{6a}
\end{equation}
For such $x$ and $\lambda$, if we assume both
\begin{equation}
\inf_{ \widehat\Omega\setminus\{0\}} \left[u- \Phi(\xi, x, \lambda;
\cdot)\right]=0 \label{6b}
\end{equation}
and
$$\liminf_{ |y|\to 0}\left[u(y)- \Phi(\xi, x, \lambda; y)\right]>0,
$$
we would have, for some $\epsilon, \epsilon'>0$,
\begin{equation}
u(y)- \Phi(\xi, x, \lambda; y)>\epsilon',\qquad \forall\ 0<|y|\le
\epsilon. \label{7a}
\end{equation}
Let
$$
\widetilde u(y):=\frac 1\lambda u(y), \qquad\widetilde
\xi(y):=\xi(x+y).
$$
We know from (\ref{6a}) and (\ref{7a}) that $ \widetilde
u>\widetilde \xi\quad\mbox{on}\
\partial(\Omega\setminus  \overline{B_\epsilon})$
i.e. $\displaystyle{ v<\lambda^{-\frac 2{n-2} }w(x+\cdot) }$ on $
\partial(\Omega\setminus  \overline{B_\epsilon})$.  Since
$\displaystyle{ \lambda^{-\frac 2{n-2} }w(x+\cdot) }$ is still a
viscosity supersolution of (\ref{L1-1}), while the singular set of
$v$ in $\Omega\setminus  \overline{B_\epsilon}$ is $S_{m-1}$ which
contains only $m-1$ points, we have, by the induction hypothesis,
$$
\inf_{   (\Omega\setminus  \overline{B_\epsilon})\setminus S_{m-1}
}\left[\lambda^{-\frac 2{n-2} }w(x+\cdot)-v\right]>0
$$
i.e.
$$
\inf_{   (\Omega\setminus  \overline{B_\epsilon})\setminus S_{m-1}
}\left[u-\Phi(\xi, x, \lambda; \cdot)\right]>0.
$$
This and (\ref{7a}) violate (\ref{6b}).  Impossible. We have proved
that (\ref{6b}) implies $$\liminf_{|y|\to 0} \left[ u(y)-\Phi(\xi,
x, \lambda; y)\right]=0 .$$  Therefore we can apply theorem 1.6 in
\cite{Li2005c} to obtain, in view of (\ref{5b}), $ u=\xi=\xi(0)$
near the origin.  Lemma \ref{lemA5} is established.

\bigskip

Because of Lemma \ref{lemA5}, \begin{equation} v=w=w(0)
\qquad\mbox{in}\ B_\epsilon \label{9a}
\end{equation}
and therefore $v$ is a viscosity solution of $\lambda(A_v)\in
\Gamma$ in $B_\epsilon$.  Thus $v\in C^{0,1}_{loc}(\overline
\Omega\setminus S_{m-1})$ is a viscosity subsolution of
$\lambda(A_v)\in \partial \Gamma$ in $\Omega\setminus S_{m-1}$.  By
the induction hypothesis, we have
$$
\inf_{ \Omega\setminus S_{m-1}}(w-v)>0.
$$
This violates (\ref{9a}).  Impossible. Step 2 is established.
  We have therefore proved Proposition \ref{propA5}.

\section{Proof of Theorem \ref{thm2}}

\noindent{\bf Proof of Theorem \ref{thm2}.}\
Suppose the contrary of (\ref{cc2}),
then in $B_{2}$,
the ball  in $\Bbb R^n$ of radius $2$ and centered at the origin,
 there exists a sequence of
 $C^4$ functions $\{u_i\}$, $C^2$ functions $\{h_i\}$ and
  $n\times n$ symmetric
positive definite $C^4$ matrix functions
$(a_{lm}^{(i)}(x))$
  satisfying,
 for some $\bar a>0$,
\begin{equation}
\frac 1{\bar a} |\xi|^2\le a_{lm}^{(i)}(x)\xi^l\xi^m\le \bar a
 |\xi|^2,
\qquad \forall\ x\in B_2, \ \xi\in \Bbb R^n,
\label{cond1}
\end{equation}
$$
\|a_{lm}^{(i)}\|_{ C^{4}(B_2) }, \|h_i\|_{ C^2(B_2) }\le \bar a,
$$
\begin{equation}
0< u_i\le \bar a\qquad \mbox{on}\ B_2,
\label{cond3}
\end{equation}
and, for the Riemannian metric
\begin{equation}
g_i:= a_{lm}^{(i)}(x)dx^ldx^m,
\label{2-1}
\end{equation}
\begin{equation}
f(\lambda(A_{ u_i^{  \frac 4{n-2} }g_i }))=h_i, \quad
\lambda(A_{ u_i^{  \frac 4{n-2} }g_i })\in \Gamma,
\qquad \mbox{in}\ B_2,
\label{c1}
\end{equation}
$$
\sup_{B_{\frac 12} }|\nabla \log u_i|\to\infty.
$$

It follows, for some $x_i\in B_1$,  that
$$
(1-|x_i|) |\nabla \log u_i(x_i)|=
\max_{ |x|\le 1} (1-|x|)|\nabla \log u_i(x)|\to \infty,
$$
where $|x|:=\sqrt{ \sum_{l=1}^n (x_l)^2 }$.
Let
$\sigma_i:= (1-|x_i|)/2$ and $
\epsilon_i:=
(2|\nabla \log u_i(x_i)|)^{-1}$.
   Then
\begin{equation}
\frac  {\sigma_i}{ \epsilon_i}\to\infty,\qquad
 2|\nabla \log u_i(x_i)|
\ge  |\nabla \log u_i(x)|,\ \ \
\forall\ |x-x_i|<\sigma_i.
\label{8-1}
\end{equation}
Consider
\begin{equation}
v_i(y):= \frac 1{u_i(x_i)}u_i(x_i+\epsilon _iy),
\qquad |y|<\frac {\sigma_i}{ \epsilon_i}.
\label{vi}
\end{equation}
Then $
v_i(0)=1$ and,
by  (\ref{8-1}) and the definition of $\epsilon_i$,
\begin{equation}
|\nabla \log v_i(y)|\le 2|\nabla \log v_i(0)|=1,
\qquad\quad \forall\ |y|<\frac {\sigma_i}{ \epsilon_i}.
\label{9-2}
\end{equation}
Thus for any $\beta>1$ there exists some
positive constant $C(\beta)$,
independent of $i$,
such that
\begin{equation}
\frac 1{C(\beta)}
\le v_i(y)\le C(\beta)\qquad
\forall\ |y|<\beta.
\label{9-3}
\end{equation}

For  $ g^{(i)}=a^{(i)}_{lm}(x_i+\epsilon_i y)dy^l dy^m$,
 $\gamma_i:= u_i(x_i)^{ - \frac 4{n-2} }\epsilon_i^{-2}
\to\infty$, and $x=x_i+\epsilon_i y$,
\begin{equation}
f(\gamma_i\lambda(A_{v_i(y)^{\frac 4{n-2} } g^{(i)} }))
=
f(\lambda(A_{u_i(x)^{\frac 4{n-2} } g_i}))
=h_i,\qquad |y|<\frac {\sigma_i}{\epsilon_i}.
\label{a3}
\end{equation}
By the proof of (1.39) in \cite{LL},
applied to
$
f(\gamma_i \cdot)$, we have, for a possibly larger $C(\beta)$,
\begin{equation}
|\nabla^2 v_i(y)|\le C(\beta)\qquad\forall\
|y|\le \beta.
\label{a4}
\end{equation}
Passing to a subsequence,
$
v_i\to v\ \mbox{in}\ C^{1,\alpha}_{loc}(\Bbb R^n)\
\mbox{for all}\ 0<\alpha<1,
$
where  $v$ is a  positive  function in $C^{1,1}_{loc}(\Bbb R^n)$
satisfying
$
|\nabla v(0)|=\frac 12.$  In particular, $v$ can not be a constant.

By
 (\ref{9-3}), (\ref{9-2}) and (\ref{a4}),
$
|\lambda( A_{v_i(y)^{\frac 4{n-2} } g^{(i)} })|
\le C(\beta)\ \forall\
|y|\le \beta.
$
This and (\ref{a3}) imply, in view of
(\ref{1-4}) and (\ref{1-5}), that
$
\lim_{i\to \infty}
f( \lambda( A_{v_i(y)^{\frac 4{n-2} } g^{(i)} }))=0.
$
Therefore
 $v$ is a $C^{1,1}_{loc}$ solution of $f(\lambda(A^v))=0$ in $\Bbb R^n$.
 By theorem 1.3 in \cite{Li2005c},
$v$ is identically a constant. A contradiction.
Theorem \ref{thm2} is established.

\medskip

It is easy to see from the proof that in Remark \ref{aabbcc}
assumption (\ref{3-3}) can be replaced by the following weaker one:
$$
\inf_{  \lambda\in \Gamma, |\lambda|\ge \frac 1\delta } \left(
|\lambda|^{2-\delta}\sum_i f_{\lambda_i}(\lambda)\right) \ge
\delta,\ \ \mbox{for some}\ \delta>0.
$$

\section{Proof of Theorem \ref{local2}}

In this section we prove Theorem \ref{local2}. We first introduce
some notations. Let  $v$ be a locally Lipschitz function in
 some open subset
 $\Omega$ of $\Bbb R^n$.
For $0<\alpha<1$,  $x\in \Omega$ and $0<\delta <dist(x,\partial\Omega)$,
let
$$
[v]_{\alpha, \delta}(x):=\sup_{ 0<|y-x|<\delta} \frac { |v(y)-v(x)|
}{ |y-x|^\alpha},
$$
$$
\delta(v,x;\Omega,\alpha):=
\left\{
\begin{array}{ll}
\infty&\mbox{if}\ [v]_{\alpha, dist(x, \partial \Omega)}(x)<1,\\
\mu\ \mbox{where} 0<\mu\le dist(x, \partial \Omega),
\mu^\alpha[v]_{\alpha, \mu}(x)=1 & \mbox{if}\ [v]_{\alpha,dist(x,
\partial \Omega)}\ge 1.
\end{array}
\right.
$$

\noindent{\bf Proof of Theorem \ref{local2}.}\ Before
establishing the gradient estimate of $\log u$, we first prove the following
H\"older estimates:
\begin{equation}
\sup_{|y|, |x|<r,  |y-x|< 2r}
\frac {  |\log u(y)-\log u(x)| }{ |y-x|^\alpha }
\le C(\alpha), \qquad \forall\ 0<\alpha <1.
\label{holder}
\end{equation}
Suppose the contrary of (\ref{holder}), then for some $0<\alpha<1$,
 there exist,  in
$B_2$,
 $C^3$ functions $\{u_i\}$, $C^1$ functions $\{h_i\}$ and
  $n\times n$ symmetric
positive definite $C^3$ matrix functions
$(a_{lm}^{(i)}(x))$,
  satisfying,
 for some $\bar a>0$, (\ref{cond1}), (\ref{cond3}) and
$$
\|a_{lm}^{(i)}\|_{ C^{3}(B_2) }, \|h_i\|_{ C^1(B_2) }\le \bar a,
$$
and (\ref{c1}) holds with
$g_i$ given by (\ref{2-1}),
but
$$
\inf_{x\in B_{\frac 12} }\delta(\log u_i, x)
\to 0,
$$
where
$$
\delta(\log u_i, x):=\delta(\log u_i, x; B_2, \alpha).
$$
It follows, for some $x_i\in B_1$,
$$
\frac{ 1-|x_i|}{  \delta(\log u_i, x_i)  }
=
\max_{|x|\le 1}
\frac{ 1-|x|}{  \delta(\log u_i, x)  }
\to \infty.
$$
Let
\begin{equation}
\sigma_i:= \frac { 1-|x_i|}{  2 },
\qquad
\epsilon_i:= \delta(\log u_i, x_i).
\label{b4}
\end{equation}
Then
\begin{equation}
\frac {\sigma_i} {   \epsilon_i  }\to \infty,
\qquad \epsilon_i\to 0,
\label{19-1}
\end{equation}
and
\begin{equation}
\epsilon_i\le 2\delta(\log u_i, z)\qquad \forall\ |z-x_i|<\sigma_i.
\label{18-1}
\end{equation}

Let $v_i$ be defined as in (\ref{vi}) with the new $\epsilon_i$
above.   By the definition of $\delta(\log u_i, x_i)$,
$$
 [\log v_i]_{\alpha, 1}(0)= \epsilon_i^\alpha[\log
u_i]_{\alpha, \epsilon_i}(x_i) =  \delta(\log u_i, x_i)^\alpha [\log
u_i]_{\alpha,  \delta(\log u_i, x_i)}(x_i)=1. $$

For any $\beta>1$ and
$|x|<\beta$, we have, in view of (\ref{19-1}), (\ref{18-1})
and the triangle inequality,
that for large $i$,
\begin{eqnarray*}
&&|\log u_i(z)-\log u_i(x_i+\epsilon_i x)|\\
&\le &|\log u_i(z)-\log u_i(\frac 12(z+x_i+\epsilon_i x))|
+|\log u_i(\frac 12(z+x_i+\epsilon_i x)) -\log u_i(x_i+\epsilon_i x)|,
\end{eqnarray*}
$$
|z-(x_i+\epsilon_i x)|=
2|z-\frac 12(z+x_i+\epsilon_i x)|=
2|\frac 12(z+x_i+\epsilon_i x) -(x_i+\epsilon_i x)|,
$$
\begin{eqnarray}
[\log v_i]_{ \alpha, 1}(x)
&=& \epsilon_i^\alpha [\log u_i]_{\alpha, \epsilon_i}(x_i+\epsilon_i x)
\nonumber\\
&\le& 2^{-\alpha}
\epsilon_i^ \alpha
\left( \sup _{ |z-(x_i+\epsilon_i x)|<\epsilon_i }
[\log u_i]_{\alpha, \frac {\epsilon_i}2 }(z)+
   [\log u_i]_{\alpha, \frac {\epsilon_i}2}(x_i+\epsilon_i x)
\right)\nonumber\\
&\le &
C(\beta)\bigg( \sup _{ |z-(x_i+\epsilon_i x)|<\epsilon_i }
\delta(\log u_i, z)^\alpha
[\log u_i]_{\alpha,   \delta(\log u_i, z)}(z)\nonumber\\
&&+
\delta(\log u_i, x_i+\epsilon_i x)^\alpha
 [\log u_i]_{\alpha, \delta(\log u_i, x_i+\epsilon_i x)}(x_i+\epsilon_i x)
\bigg)\le C(\beta).\nonumber
\end{eqnarray}
This implies (\ref{9-3}) for any $\beta>1$. By Theorem D, we have,
for any $\beta>1$,
$$
|\nabla v_i(y)|\le C(\beta)\qquad \forall\ |y|<\beta.
$$
Passing to a subsequence,
$$
v_i\to v\qquad \mbox{in}\ C^{\gamma}_{loc}(\Bbb R^n)\ \mbox{for
all}\ \alpha<\gamma<1,
$$
where  $v$ is a  positive  function in $C^{0,1}_{loc}(\Bbb R^n)$
satisfying
$
[\log v]_{\alpha, 1}(0)=1.
$  In particular, $v$ can not be a constant.

Clearly,  (\ref{a3}) holds
with  the  new $\epsilon_i$ given in
(\ref{b4}).
Thus, by (\ref{1-5new}) and
(\ref{a3}),
$$
\lim_{i\to \infty}
f\left( \lambda\left( A_{v_i(y)^{\frac 4{n-2} } g^{(i)} }\right)\right)=0.
$$
It is easy to see that $w:= v^{ -\frac 2{n-2}}$
 is a positive locally Lipschitz viscosity solution of
$\lambda(A_w)\in \partial \Gamma$ in $\Bbb R^n$.  By Theorem
\ref{thm-vis}, $v$ is a constant. A contradiction. The H\"older
estimate (\ref{holder}) is established.

Now we establish the gradient estimate (\ref{cc2})
based on the Holder estimates.
The H\"older estimate (\ref{holder}) yields
the Harnack inequality:
$$
\sup_{ B_{2r}}u
\le C \inf_{ B_{2r} }u.
$$
Consider
$$
w:= \frac 1{u(0)}u.
$$
The equation of $w$ on $B_{3r}$  is
$$
f(\lambda(A_{ w^{  \frac 4{n-2} }g }))=u(0)^{ \frac 4{n-2} }h, \quad
\lambda(A_{ w^{  \frac 4{n-2} }g })\in \Gamma,
$$
and $w$ satisfies
$$
\frac 1C \le w\le C\qquad \mbox{in}\ B_{2r}.
$$
Since $u(0)$ is bounded from above, we have, using Theorem D, that
$$
|\nabla u|\le C\quad \mbox{in}\ B_r.
$$
Theorem \ref{local2} is established.

\section{Proof of Theorem \ref{thm16} and Theorem \ref{boundary}}

\begin{prop} Let $u^+\in C^{0,1}(\overline{B_1^+})$ and
 $u^-\in C^{0,1}(\overline{B_1^-})$ be two positive function
 satisfying $u^+=u^-$ on $\partial'B_1^+$.
We assume that \begin{equation} u^+ \ \mbox{is a viscosity
supersolution (subsolution) of}\
 \lambda(A^{u^+})\in
\partial\Gamma\ \mbox{in}\ B_1^+
\label{aaa}
\end{equation}
\begin{equation}
\frac{\partial u^+}{\partial x_n}\le \ (\ge)\ 0 \ \mbox{on}\
\partial'B_1^+\
\mbox{in the viscosity sense}, \label{15anew}
\end{equation}
\begin{equation} u^- \ \mbox{is a viscosity
supersolution (subsolution) of}\
 \lambda(A^{u^-})\in
\partial\Gamma\ \mbox{in}\ B_1^+
\label{aaanew}
\end{equation}
\begin{equation}
\frac{\partial u^-}{\partial x_n}\ge\ (\le)\ 0 \ \mbox{on}\
\partial'B_1^-\
\mbox{in the viscosity sense}. \label{15anewnew}
\end{equation}
Then
$$
\widetilde u(x', x_n):= \left\{
\begin{array}{rl}
u^+(x', x_n)& \mbox{if}\ x_n\ge 0,\\
u^-(x', x_n)& \mbox{if}\ x_n<0 \end{array} \right.
$$
is a $C^{0,1}$ viscosity supersolution (subsolution) of
$\lambda(A^{\widetilde u})\in
\partial\Gamma$ in $B_1$. \label{LemmaA}
\end{prop}

A consequence of Proposition \ref{propA5} and Proposition
\ref{LemmaA} is

\begin{cor}
Let $\Omega^+\subset \Bbb R^n_+$ be a bounded open set. For $m$
points $S_m:= \{P_1, \cdots, P_m\}\subset \Omega^+\cup
\partial'\Omega^+$, $m\ge 0$, let $u\in
C^{0,1}(\overline{\Omega^+}\setminus S_m)$ and $v\in
C^{0,1}(\overline{\Omega^+})$ be positive functions. Assume that
$$
u \ \mbox{is a viscosity supersolution of}\ \lambda(A^u)\in \partial
\Gamma\ \mbox{in}\ \Omega^+\setminus S_m, $$
$$
 v \ \mbox{is a viscosity subsolution of}\
\lambda(A^v)\in \partial \Gamma\ \mbox{in}\ \Omega^+, $$
$$\frac{\partial u}{\partial x_n}\le 0 \le \frac{\partial v}{\partial
x_n}\ \mbox{on}\
\partial'\Omega^+\
\mbox{in the viscosity sense},
$$
$$
 u>v\qquad\mbox{on}\
\partial''\Omega^+.
$$
Then
\begin{equation}
\inf_{\overline{\Omega^+}\setminus S_m}(u-v)>0. \label{15c}
\end{equation}
\label{prop14}
\end{cor}

\noindent{\bf Proof of Corollary \ref{prop14} using Proposition
\ref{propA5} and Proposition \ref{LemmaA}.}\ Let $u^:=u$, $u^-(x',
x_n):=u(x', -x_n)$, $v^+:=v$, $v^-(x', x_n):=v(x', -x_n)$,

$$
\widetilde u(x', x_n):= \left\{
\begin{array}{rl}
u^+(x', x_n)& \mbox{if}\ x_n\ge 0,\\
u^-(x', x_n)& \mbox{if}\ x_n<0 \end{array} \right.
$$
and

$$
\widetilde v(x', x_n):= \left\{
\begin{array}{rl}
v^+(x', x_n)& \mbox{if}\ x_n\ge 0,\\
v^-(x', x_n)& \mbox{if}\ x_n<0 \end{array} \right.
$$
An application of Proposition \ref{propA5} to $v^{-\frac 2{n-2}}$
and $u^{-\frac 2{n-2}}$, in view of Proposition \ref{LemmaA} and
Remark \ref{rem1.10}, yields (\ref{15c}).

\medskip

\noindent{\bf Proof of Proposition \ref{LemmaA}.} We first prove the
proposition under, instead of (\ref{15anew}) and (\ref{15anewnew}),
\begin{equation}
\frac{\partial u^+}{\partial x_n}<0 \ \mbox{on}\
\partial'B_1^+\ \mbox{and}\ \frac{\partial u^-}{\partial x_n}>0 \ \mbox{on}\
\partial'B_1^-\
\mbox{in the viscosity sense}. \label{26b}
\end{equation}
Let $\bar x\in B_1$, $\psi\in C^2(B_1)$, $\widetilde u(\bar
x)=\psi(\bar x)$ and, for some $0<\delta<1-|\bar x|$, $u(x)\ge
\psi(x)$ for all $|x-\bar x|<\delta$.  We need to show that
$$\lambda(A^\psi(\bar x))\in \Bbb R^n\setminus \Gamma.
$$
If $\bar x$ does not belong to $\partial' B_1^+$, this is obvious
because of (\ref{aaa}) and (\ref{aaanew}).  So we only need to show
that $\bar x$ does not belong to $\partial' B_1^+$.  Indeed, if
$\bar x\in
\partial' B_1^+$, then, since $\bar x_n=0$,
$$
u^+(\bar x)=u^-(\bar x)=\psi(\bar x),\  u^+, u^-\ge \psi\
\mbox{near}\ \bar x.
$$
Thus, by (\ref{26b}),
$$
\frac{\partial \psi}{\partial x_n}(\bar x)<0, \mbox{and}\
\frac{\partial \psi}{
\partial x_n}(\bar x)>0.
$$
A contradiction.

Now we prove the proposition under (\ref{15anew}) and
(\ref{15anewnew}). We will only give the proof when $u^+$ and $u^-$
are viscosity supersolutions, since the proof is essentially the
same when they are subsolutions.  We start with a first variation of
the operator $A^u$ together with the Neumann boundary condition.

\begin{lem}  Let $\Omega^+\subset \Bbb R^n_+$ be a bounded open set,
$w\in C^2(\Omega^+)\cap C^1(\Omega^+\cup \partial'\Omega^+)$
satisfy, for some constant $c_1>0$,
$$
w\ge c_1\ \ \mbox{in}\ \Omega^+
$$
and let
$$\varphi^{\pm}(x):=e^{ \delta |x|^2 \pm \delta^2 x_n}.
$$
Then there exists some constant $\delta>0$, depending only on
$\sup\{ |x|\ |\ x\in \Omega^+\}$, and there exists $\bar
\epsilon>0$, depending only on $\delta, c_1$ and $\sup \{ |x|\ |\
x\in \Omega^+\}$, such that for any $0<\epsilon<\bar \epsilon$,
$$
A_{w+\epsilon\varphi^{\pm}}\ge \left(1+\epsilon\frac{\varphi^{\pm}}
w\right)A_w+\frac{\epsilon\delta}2 \varphi^{\pm}wI\qquad \mbox{in}\
\Omega^+,
$$
 $$
A_{w-\epsilon\varphi^{\pm}}\le \left(1-\epsilon\frac{\varphi^{\pm}}
w\right)A_w-\frac{\epsilon\delta}2 \varphi^{\pm}wI\qquad \mbox{in}\
\Omega^+,
$$
$$
\frac{\partial}{\partial x_n}(w+\epsilon\varphi^{\pm})=
\frac{\partial w}{\partial x_n}\pm \epsilon\delta\qquad\mbox{on}\
\partial'\Omega^+,
$$
$$
\frac{\partial}{\partial x_n}(w-\epsilon\varphi^{\pm})=
\frac{\partial w}{\partial x_n}\mp \epsilon\delta\qquad\mbox{on}\
\partial'\Omega^+,
$$
\label{lem29}
\end{lem}

\noindent{\bf Proof.}\  It is very similar to that of lemma 3.7 in
\cite{Li2005c}, we omit the details.

\bigskip

Let $u^+$ be the supersolution in Proposition \ref{LemmaA}, set
$$
\xi^+:= (u^+)^{-\frac 2{n-2}}, \quad
\xi_\epsilon^+:=\xi^++\epsilon\varphi^+, \quad
u_\epsilon^+:=(\xi_\epsilon^+)^{-\frac {n-2}2 }.
$$
We will prove that
\begin{equation}
u_\epsilon^+\ \mbox{is a viscosity supersolution of}\
\lambda(A^{u_\epsilon^+})\in
\partial \Gamma\ \mbox{in}\ B_1^+,
\label{31a}
\end{equation}
and
\begin{equation}
\frac{\partial u_\epsilon^+}{\partial x_n}<0\ \mbox{on}\
\partial' B_1^+\ \mbox{in the viscosity sense}.
\label{31b}
\end{equation}

Let $\bar x\in B_1^+$, $\psi\in C^2(B_1^+)$, $u_\epsilon^+(\bar
x)=\psi(\bar x)$ and $u_\epsilon\ge \psi$ near $\bar x$.  Then, with
$\eta:=\psi^{ -\frac 2{n-2} }$,
$$
\xi^+=\eta-\epsilon\varphi^+\ \mbox{at}\ \bar x\ \mbox{and}\
\xi^+\le \eta-\epsilon\varphi^+\ \mbox{near}\ \bar x.
$$
By Remark \ref{rem1.10}, $\xi^+$ is a viscosity subsolution of
$\lambda(A_{\xi^+})\in
\partial \Gamma$, and therefore
$$
\lambda(A_{\eta-\epsilon\varphi^+}(\bar x))\in \overline \Gamma.
$$
By Lemma \ref{lem29},
$$
A_{\eta-\epsilon\varphi^+}(\bar x)< \left(1-\epsilon
\frac{\varphi^+}\eta\right)(\bar x)A_\eta(\bar x)
$$
which implies, for small $\epsilon$,
$$
A^\psi(\bar x)=A_\eta(\bar x)\in \Gamma.
$$
We have proved (\ref{31a}).  To prove (\ref{31b}), let $\bar x\in
\partial' B_1^+$, $\psi\in C^1(\overline {B_1^+})$,
$u_\epsilon^+(\bar x)=\psi(\bar x)$ and $u_\epsilon\ge \psi$ near
$\bar x$.  It follows that
$$
u^+=\left[ \psi^{-\frac 2{n-2} }-\epsilon \varphi^+\right]^{-\frac
{n-2}2}\ \ \mbox{at}\ \bar x \ \ \mbox{and}\ u^+\ge \left[
\psi^{-\frac 2{n-2} }-\epsilon \varphi^+\right]^{-\frac {n-2}2}\ \
\mbox{near}\ \bar x.
$$
Since $\frac{\partial u^+}{\partial x_n}\le 0$ on $\partial' B_1^+$
in the viscosity sense, we have
$$
0\ge  \frac{\partial}{\partial x_n} \left[ \psi^{-\frac 2{n-2} }
-\epsilon \varphi^+\right]^{-\frac {n-2}2} \bigg|_{ x=\bar x} =
\frac{\partial \psi}{\partial x_n}\left[ 1+O(\epsilon)\right]
+\frac{n-2}2 \epsilon \delta^2 \psi^{\frac n{n-2}}+O(\epsilon^2).
$$
So, for small $\epsilon$, we have $\frac{\partial \psi}{\partial
x_n}(\bar x)<0$.  We have proved (\ref{31b}).

Similarly we set for $u^-$
$$
\xi^-= (u^-)^{-\frac 2{n-2}}, \quad
\xi_\epsilon^-:=\xi^-+\epsilon\varphi^-, \quad
u_\epsilon^-:=(\xi_\epsilon^-)^{-\frac {n-2}2 },
$$
and can prove
$$
u_\epsilon^-\ \mbox{is a viscosity supersolution of}\
\lambda(A^{u_\epsilon^-})\in
\partial \Gamma\ \mbox{in}\ B_1^-,
$$ and
\begin{equation}
\frac{\partial u_\epsilon^-}{\partial x_n}>0\ \mbox{on}\
\partial' B_1^-\ \mbox{in the viscosity sense}.
\label{31bnew} \end{equation}
 It is clear that
$u_\epsilon^+=u_\epsilon^-$ on $\partial B_1^+$.

Since we now have the strict inequalities (\ref{31b}) and
(\ref{31bnew}),
$$
\widetilde u_\epsilon(x', x_n):= \left\{
\begin{array}{rl}
u_\epsilon^+(x', x_n)& \mbox{if}\ x_n\ge 0,\\
u_\epsilon^-(x', x_n)& \mbox{if}\ x_n\le0 \end{array} \right. $$
 is
a viscosity supersolutions of $\lambda(A^{\widetilde u_\epsilon})\in
\partial\Gamma$ in $B_1$. Since $\widetilde u_\epsilon\to \widetilde
u$ in $C^0_{loc}(B_1)$, we have, by standard arguments, $\widetilde
u$ is a viscosity supersolution of $\lambda(A^{\widetilde u})\in
\partial \Gamma$ in $B_1$.
Proposition \ref{LemmaA} is established.

\noindent{\bf Proof of Theorem \ref{thm16}.}\  By Proposition
\ref{LemmaA},
$$
\widetilde u(x', x_n):= \left\{
\begin{array}{rl}
u(x', x_n)& \mbox{if}\ x_n\ge 0,\\
u(x', -x_n)& \mbox{if}\ x_n\le0 \end{array} \right. $$ satisfies the
hypothesis of Theorem \ref{thm-vis} and therefore is a constant.

\noindent{\bf Proof of Theorem \ref{boundary}.}\  The proof is
similar to that of Theorem \ref{local2}.  Let $O_3$ be an open set
of $M$ satisfying $\overline O_3\subset O_3\subset \overline
O_3\subset O_1$.  We first establish
\begin{equation}
\label{holder1} \sup_{y, x\in O_3,  dist(y,x)< 2r} \frac {  |\log
u(y)-\log u(x)| }{ dist(y,x)^\alpha } \le C(\alpha), \qquad \forall\
0<\alpha <1.
\end{equation}
Suppose the contrary of (\ref{holder1}), then for some $0<\alpha<1$,
 there exist,  in
$\overline B_2^+\subset \Bbb R^n$,
 $C^3$ functions $\{u_i\}$, $C^1$ functions $\{\psi_i\}$ and
 $\{\eta_i\}$, and
  $n\times n$ symmetric
positive definite $C^3$ matrix functions $(a_{lm}^{(i)}(x))$,
  satisfying,
 for some $\bar a>0$, (\ref{cond1}) and (\ref{cond3})
 in $\overline B_2^+$,  and
$$
\|a_{lm}^{(i)}\|_{ C^{3}(B_2^+) }, \|\psi_i\|_{ C^1(B_2^+) },
\|\eta_i\|_{ C^1(B_2^+) }\le \bar a,
$$
$$
\left\{
\begin{array}{ll}
f(\lambda(A_{u_i^{\frac 4{n-2}}g_i}))=\psi_i, &
\lambda(A_{u_i^{\frac
4{n-2}}g_i})\in \Gamma\ \mbox{on}\ B_2^+,\\
-\frac{ \partial u_i}{\partial \nu_{g_i}}+\frac {n-2}2
h_{g_i}u_i=\eta_iu_i^{ \frac n{n-2}},& \mbox{on}\ \partial' B_2^+,
\end{array}
\right. $$
where $g_i$ is given by (\ref{2-1}), but
$$
\inf_{x\in B_{\frac 12}^+ }\delta(\log u_i, x) \to 0,
$$
where
$$
\delta(\log u_i, x):=\delta(\log u_i, x; B_2^+, \alpha).
$$
It follows, for some $x_i\in B_1^+\cup \partial' B_1^+$,
$$
\frac{ 1-|x_i|}{  \delta(\log u_i, x_i)  } = \max_{x\in \overline
{B_1^+}} \frac{ 1-|x|}{ \delta(\log u_i, x)  } \to \infty.
$$
Let $ \sigma_i$ and $ \epsilon_i$ be defined as in (\ref{b4}). Then
they satisfy (\ref{19-1}) and
$$
\epsilon_i\le 2\delta(\log u_i, z)\qquad \forall\ z\in
B_{\sigma_i}(x_i)\cap B_1^+. $$ Let
$$
v_i(y):= \frac 1{u_i(x_i)}u_i(x_i+\epsilon _iy), \qquad |y|<\frac
{\sigma_i}{ \epsilon_i}, y_n>-T_i:= -\frac 1{\epsilon_i} (x_i)_n.
$$
After passing to a subsequence, either $\lim_{i\to
\infty}(-T_i)=-T>-\infty$ or $\lim_{i\to \infty} (-T_i)=-\infty$.
Following, with obvious modification, the arguments in the proof of
Theorem \ref{local2}, we see, passing to another subsequence, that
either
$$ v_i(\cdot +(0', -T_i))\to v\ \ \mbox{in}\
C^\gamma_{loc}(\overline {\Bbb R^n_+})\ \ \mbox{for all}\
0<\gamma<1,
$$
for some positive locally Lipschitz viscosity solution $v$ of
(\ref{16a}) and (\ref{16b}) satisfying $[\log v]_{ \alpha, 1}(0',
T)=1$, or
$$
v_i\to v\ \ \mbox{in}\ C^\gamma_{loc}(R^n)\ \ \mbox{for all}\
0<\gamma<1,
$$
for some positive locally Lipschitz viscosity solution of
$\lambda(A^v)\in \partial \Gamma$ in $\Bbb R^n$ satisfying $[\log
v]_{ \alpha, 1}(0)=1$.   By our Liouville theorems $v$ must be a
constant.  But $[\log v]_{ \alpha, 1}(0', T)=1$ or $[\log v]_{
\alpha, 1}(0)=1$  does not allow $v$ to be a constant.  A
contradiction. Theorem \ref{boundary} is established.

\section{Appendix A}

We give in this appendix the proof of Theorem D in \cite{LL3}.  For
simplicity we present the proof on locally conformally flat
manifolds.  Namely we give another proof of Theorem C, which can
easily be extended to general Riemannian manifolds.

\noindent{\bf Another proof of Theorem C.}\  We write $v=-\frac
2{n-2}\log u$.  Then $v$ satisfies, with $\alpha=-\frac 2{n-2}\log
b$ and $\beta=-\frac 2{n-2}\log a$,
\begin{equation}
f(\lambda(W))=h,\ \ \lambda (W)\in\Gamma,\quad\mbox{in}~B_3,
\label{citee1}
\end{equation}
and
$$
\alpha\le v\le \beta\qquad\mbox{on}\ B_3,
$$
where
\[
W:=(W_{ij})=e^{2v}(v_{ij}+v_iv_j-\frac{|\nabla v|^2}{2}\delta_{ij}).
\]
We only need to prove, for some constant $C$ depending  on $\alpha$,
$\beta$ and  $(f,\Gamma),$, that
\begin{equation}
|\nabla v|\le C\qquad\mbox{on}~B_1. \label{ZZZ}
\end{equation}

Fixing some small constants $\epsilon, c_1>0$, depending only on
$\alpha$ and $\beta$, such that the function $\phi(s):=\epsilon
e^{-2s}$ satisfies
\begin{equation}
-\frac 12\phi'\ge c_1,\quad \phi''+\phi'-(\phi')^2\ge 0,\quad\mbox
{on}~[\alpha, \beta]. \label{claim}
\end{equation}

Let $\rho\ge 0$ be a smooth function taking value $1$ in $B_1$ and
$0$ outside $B_2$. It is known that $\rho$ satisfies $|\nabla
\rho|^2\le C_1$.
 Consider
$$
G=\rho e^{\phi(v)}|\nabla v|^2. $$ Estimate (\ref{ZZZ}) is
established if we can show that $G\le C$ on $\bar B_2$. Let
$G(x_0)=\max\limits_{\bar B_2}G$ for some $x_0\in\bar B_2$. Clearly
$x_0\in B_2$. After a rotation of the axis if necessary, we may
assume that $W(x_0)$ is a diagonal matrix. In the following, we use
subscripts of a function to denote derivatives. For example,
$G_i=\partial_{x_i}G$, $G_{ij}=\partial_{x_ix_j}G$, and so on.
 We also use the notation $f^i:=\frac{\partial f}{\partial
\lambda_i}$.

Applying $\partial_{x_k}$ to (\ref{citee1}) leads to
\begin{equation}
f^i W_{iik}=0. \label{diff} \end{equation} By calculation,
\begin{eqnarray*}
G_i&=&2\rho e^{\phi}v_{ki}v_k+\rho\phi'e^{\phi}|\nabla v|^2
v_i+e^{\phi}|\nabla v|^2\rho_i\\
&=&2\rho e^{\phi}v_{ki}v_k+(\phi' v_i+\frac{\rho_i}{\rho})G\\
\end{eqnarray*}
At $x_0$, we have $G_i=0$. Equivalently, we have
\begin{equation}
2v_{ki}v_k=-\phi' |\nabla v|^2 v_i-\frac{\rho_i}{\rho}|\nabla
v|^2,\quad\forall~ 1\le i\le n. \label {G_i}
\end{equation}
Take the second covariant derivative of $G$ and evaluate at $x_0$,
\begin{eqnarray*}
0\ge (G_{ij})&=&2v_{kij}v_k
e^{\phi}\rho+2v_{ki}v_{kj}e^{\phi}\rho+2v_{ki}v_ke^{\phi}\phi'
v_j\rho
+2v_{ki}v_k e^{\phi}\rho_j\\
&&+(\phi''v_iv_j+\phi'v_{ij}+\frac{\rho\rho_{ij}-\rho_i\rho_j}{\rho^2})\rho
e^{\phi}|\nabla v|^2,
\end{eqnarray*}
therefore, at $x_0$,
\begin{eqnarray}
0&\ge& e^{-\phi}f^iG_{ii} \nonumber \\
&=&2\rho f^i v_{iik}v_k+2\rho f^i
v_{ki}^2 +2\rho \phi'f^iv_{ki}v_kv_i +2f^i v_{ki}v_k \rho_i \nonumber \\
&&+f^i
(\phi''v_i^2+\phi'v_{ii}+\frac{\rho\rho_{ii}-\rho_i^2}{\rho^2})\rho
|\nabla v|^2 \nonumber \\
&=&2\rho f^i v_k \{e^{-2v}W_{ii}-v_i^2+\frac{|\nabla
v|^2}{2}\delta_{ii}\}_k+2\rho f^i v_{ki}^2
\nonumber \\
&&-\rho\phi' f^i(|\nabla v|^2\phi'v_i^2+\frac {|\nabla v|^2\rho_i
v_i}{\rho}) -f^i (|\nabla v|^2\phi'v_i+\frac {|\nabla
v|^2\rho_i}{\rho})\rho_i
+\rho\phi''|\nabla v|^2f^iv_i^2 \nonumber \\
&&+\rho\phi'|\nabla v|^2 f^i(e^{-2v}W_{ii}-v_i^2 +\frac{|\nabla
v|^2}{2}\delta_{ii})
+|\nabla v|^2 f^i\frac{\rho\rho_{ii}-\rho_i^2}{\rho} \nonumber \\
&=&2\rho f^i\{e^{-2v}W_{iik}v_k-2e^{-2v}|\nabla v|^2W_{ii} +\phi'
|\nabla v|^2v_i^2+\frac {\rho_i v_i}{\rho}|\nabla v|^2 -\frac 12
|\nabla v|^4\phi' \nonumber \\
&&-\frac 12 |\nabla v|^2\frac{v_k\rho_k}{\rho}
\}+2\rho f^i v_{ki}^2 \nonumber \\
&&-\rho\phi'f^i(|\nabla v|^2\phi'v_i^2+\frac {|\nabla v|^2\rho_i
v_i}{\rho}) -f^i(|\nabla v|^2\phi'v_i+\frac {|\nabla
v|^2\rho_i}{\rho})\rho_i
+\rho\phi''|\nabla v|^2f^iv_i^2 \nonumber \\
&&+\rho\phi'|\nabla v|^2 f^i(e^{-2v}W_{ii}-v_i^2 +\frac{|\nabla
v|^2}{2}\delta_{ii})+|\nabla v|^2 f^i\frac{\rho\rho_{ii}-\rho_i^2}{\rho} \nonumber \\
&=&\{-4\rho e^{-2v}|\nabla v|^2f-|\nabla v|^2 v_k\rho_k\sum\limits_i
f^i
 \nonumber \\
&&-2\phi'|\nabla v|^2f^i\rho_i v_i +\rho e^{-2v}\phi'|\nabla v|^2f \nonumber \\
&&+|\nabla v|^2f^i\frac{\rho\rho_{ii}-2\rho_i^2}{\rho}
+2|\nabla v|^2f^i\rho_i v_i\}+2\rho f^iv_{ki}^2 \nonumber \\
&&-\frac 12 \rho\phi'|\nabla v|^4\sum\limits_i f^i
+(\phi''+\phi'-(\phi')^2)\rho |\nabla v|^2 f^iv_i^2 \nonumber \\
\label {Gij}
\end{eqnarray}

In the following, we use $C_2$ to denote some positive constant
depending only on $\alpha$, $\beta$ and $(f,\Gamma)$ which may vary
from line to line. By (\ref{claim}), we derive from (\ref{Gij}) that
\begin{eqnarray*}
0\ge e^{-\phi}f^iG_{ii} &\ge & (-C_2\sqrt
{\rho}|\nabla v|^3-C_2|\nabla v|^2 +2c_1\rho |\nabla v|^4 )\sum\limits_i f^i\\
&= &|\nabla v|^2(-C_2\sqrt {\rho |\nabla v|^2} -C_2+c_1 \rho |\nabla
v|^2)\sum\limits_i f^i,
\end{eqnarray*}
which implies $\rho|\nabla v|^2(x_0)\le C_2$, so is $G(x_0)$.
Estimate (\ref{ZZZ}) is established.

\end{document}